\documentclass[10pt]{amsproc} 
\usepackage[all]{xy}
\SelectTips{eu}{}
\usepackage{ifthen} 
\usepackage{amssymb}
\calclayout
\makeatletter 
\def\serieslogo@{} 
\def\@setcopyright{} 
\makeatother
\title{Derived categories, resolutions, and Brown representability} 
\author[Henning Krause]{Henning Krause}
\address{Henning Krause\\ Institut f\"ur Mathematik\\
Universit\"at Paderborn\\ 33095 Paderborn\\ Germany.}
\email{hkrause@math.upb.de}

\newtheorem*{lem}{Lemma}

\newtheorem*{prop}{Proposition}

\newtheorem*{thm}{Theorem}

\theoremstyle{remark}

\theoremstyle{definition}
\newtheorem*{exm}{Example}

\numberwithin{equation}{section}
\newtheorem*{rem}{Remark}

\newcommand{\smatrix}[1]{\left[\begin{smallmatrix}#1\end{smallmatrix}\right]}

\renewcommand{\mod}{\operatorname{mod}\nolimits}

\newcommand{\proj}{\operatorname{proj}\nolimits}

\newcommand{\Mod}{\operatorname{Mod}\nolimits}
\newcommand{\Gr}{\operatorname{Gr}\nolimits}

\newcommand{\Hom}{\operatorname{Hom}\nolimits}
\newcommand{\RHom}{\operatorname{{\bfR}Hom}\nolimits}

\newcommand{\END}{\operatorname{\mathcal E\!\!\:\mathit n\mathit d}\nolimits}
\newcommand{\HOM}{\operatorname{\mathcal H\!\!\:\mathit o\mathit m}\nolimits}
\renewcommand{\Im}{\operatorname{Im}\nolimits}
\newcommand{\Ker}{\operatorname{Ker}\nolimits}
\newcommand{\Coker}{\operatorname{Coker}\nolimits}

\newcommand{\Ext}{\operatorname{Ext}\nolimits}

\newcommand{\Lex}{\operatorname{Lex}\nolimits}

\newcommand{\Inj}{\operatorname{Inj}\nolimits}
\newcommand{\Ab}{\mathrm{Ab}} 
\newcommand{\op}{\mathrm{op}}
\newcommand{\Id}{\mathrm{Id}} 
\newcommand{\id}{\mathrm{id}}
\newcommand{\dg}{\mathrm{dg}} 
\newcommand{\idg}{\mathrm{idg}} 
\newcommand{\pdg}{\mathrm{pdg}} 

\newcommand{\inj}{\mathrm{inj}}
\newcommand{\inc}{\mathrm{inc}}
\newcommand{\can}{\mathrm{can}}
\newcommand{\comp}{\mathop{\raisebox{+.3ex}{\hbox{$\scriptstyle\circ$}}}}
\newcommand{\lto}[1][{}]{\stackrel{#1}{\longrightarrow}} 
\renewcommand{\to}[1][{}]{\stackrel{#1}{\rightarrow}} 
\newcommand{\xto}{\xrightarrow}

\def\a{\alpha}
\def\b{\beta}
\def\e{\varepsilon}
\def\d{\delta}
\def\g{\gamma}
\def\p{\phi}
\def\r{\rho}
\def\s{\sigma}
\def\t{\tau}

\def\k{\kappa}
\def\la{\lambda}

\def\Ga{\Gamma}
\def\La{\Lambda}
\def\Si{\Sigma}

\def\A{{\mathcal A}}
\def\B{{\mathcal B}}
\def\C{{\mathcal C}}
\def\D{{\mathcal D}}

\def\I{{\mathcal I}}

\def\P{{\mathcal P}}

\def\S{{\mathcal S}}

\def\Y{{\mathcal Y}}

\def\T{{\mathcal T}}
\def\U{{\mathcal U}}
\def\V{{\mathcal V}}

\def\bbQ{\mathbb Q}

\def\bbZ{\mathbb Z}

\def\bfC{\mathbf C}
\def\bfD{\mathbf D}
\def\bfi{\mathbf i}

\def\bfK{\mathbf K}
\def\bfL{\mathbf L}
\def\bfp{\mathbf p}
\def\bfP{\mathbf P}
\def\bfR{\mathbf R}
\def\bfS{\mathbf S}
\begin{document}

\begin{abstract}
These notes are based on a series of five lectures given during the
summer school ``Interactions between Homotopy Theory and Algebra''
held at the University of Chicago in 2004.
\end{abstract}
\maketitle 

\setcounter{tocdepth}{2}
\tableofcontents
\def\addcontentsline#1#2#3{%
\addtocontents{#1}{\protect\contentsline{#2}{#3}{}}}

\section*{Introduction}

Derived categories were introduced in the 1960s by Grothendieck and
his school. A derived category is defined for any abelian category and
the idea is quite simple.  One takes as objects all complexes, and the
usual maps between complexes are modified by inverting all maps which
induce an isomorphism in cohomology.

This definition is easily stated but a better description of the maps
is needed. To this end some extra structure is introduced. A derived
category carries a triangulated structure which complements the
abelian structure of the underlying category of complexes.

The next step is to replace a complex by an injective or projective
resolution, assuming that the underlying abelian category provides
enough injective or projective objects. The construction of such
resolutions requires some machinery. Here, we use the Brown
representability theorem which characterizes the representable
functors on a triangulated category.  The proof of the Brown
representability theorem is based on an embedding of a triangulated
category into some abelian category. This illustrates the interplay
between abelian and triangulated structures.

The classical context for doing homological algebra is the module
category over some associative algebra. Here, we work more generally
over differential graded algebras or, even more generally, over
differential graded categories. Again, we prove that complexes in such
module categories can be replaced by projective or injective resolutions.

There is a good reason for studying the derived category of a
differential graded category.  It turns out that every triangulated
category which arises from algebraic constructions embeds into such a
derived category. This leads to the notion of an algebraic
triangulated category which can be defined axiomatically.

\medskip
These notes are based on a series of five lectures given during a
summer school in Chicago in 2004. The prerequisites for these lectures
are quite modest, consisting of the basic notions from homological
algebra and some experience with examples. The latter should
compensate for the fact that specific examples are left out in order
to keep the exposition concise and self-contained.

\medskip
Students and colleagues in Paderborn helped with numerous questions
and suggestions to prepare these notes; it is a pleasure to thank all
of them. Also, Apostolos Beligiannis and Srikanth Iyengar provided
many helpful comments. In addition, I am grateful to the organizers of
the summer school: Lucho Avramov, Dan Christensen, Bill Dwyer, Mike
Mandell, and most prominently, Brooke Shipley. Last but not least, I
would like to thank the participants of the summer school for their
interest and enthusiasm.

\section{Derived categories}

The derived category $\bfD(\A)$ of an abelian category $\A$ provides a
framework for studying the homological properties of $\A$. The main idea
is to replace objects in $\A$ by complexes, and to invert maps between
complexes if they induce an isomorphism in cohomology. The actual
construction of the derived category proceeds in several steps which
is reflected by the following sequence of functors.
$$\A\lto\bfC(\A)\lto\bfK(\A)\lto\bfD(\A)$$

\subsection{Additive and abelian categories}
A category $\A$ is {\em additive} if every finite family of objects
has a product, each set $\Hom_\A(A,B)$ is an abelian group, and the
composition maps
$$\Hom_\A(A,B)\times\Hom_\A(B,C)\lto\Hom_\A(A,C)$$ sending a pair
$(\phi,\psi)$ to $\psi\comp\phi$ are bilinear.

An additive category $\A$ is {\em abelian}, if every map $\p\colon
A\to B$ has a kernel and a cokernel, and if the canonical factorization
$$\xymatrix{\Ker\p\ar[r]^-{\p'}&A\ar[r]^-\p\ar[d]&
B\ar[r]^-{\p''}&\Coker\p\\ 
&\Coker\p'\ar[r]^-{\bar\p}&\Ker\p''\ar[u]}$$
of $\p$ induces  an isomorphism $\bar\p$.

\begin{exm} 
The category $\Mod\La$ of right modules over an associative ring $\La$
is an abelian category.
\end{exm}

\subsection{Categories of complexes}

Let $\A$ be an additive category. A {\em complex} in $\A$ is a
sequence of maps
$$\cdots \lto X^{n-1}\lto[d^{n-1} ]X^n\lto[d^n]X^{n+1}\lto\cdots$$
such that $d^{n}\comp d^{n-1}=0$ for all $n\in\bbZ$. We denote by
$\bfC(\A)$ the category of complexes, where a map $\p\colon X\to Y$
between complexes consists of maps $\p^n\colon X^n\to Y^n$ with
$d_Y^{n}\comp\p^n=\p^{n+1}\comp d_X^n$ for all $n\in\bbZ$. 

A map $\p\colon X\to Y$ is {\em null-homotopic} if there are maps
$\r^n\colon X^n\to Y^{n-1}$ such that $\p^n=d_Y^{n-1}\comp
\r^{n}+\r^{n+1}\comp d_X^n$ for all $n\in\bbZ$. The null-homotopic
maps form an {\em ideal} $\I$ in $\bfC(\A)$, that is, for each pair
$X,Y$ of complexes a subgroup
$$\I(X,Y)\subseteq\Hom_{\bfC(\A)}(X,Y)$$
such that any composition
$\psi\comp\p$ of maps in $\bfC(\A)$ belongs to $\I$ if $\p$ or $\psi$
belongs to $\I$.  The {\em homotopy category} $\bfK(\A)$ is the
quotient of $\bfC(\A)$ with respect to this ideal. Thus
$$\Hom_{\bfK(\A)}(X,Y)=\Hom_{\bfC(\A)}(X,Y)/\I(X,Y)$$ for every pair
of complexes $X,Y$.

Now suppose that $\A$ is abelian. Then one defines for a complex $X$
and each $n\in\bbZ$ the {\em cohomology} $$H^nX=\Ker d^n/\Im
d^{n-1}.$$ A map $\p\colon X\to Y$ between complexes induces a map
$H^n\p\colon H^nX\to H^nY$ in each degree, and $\p$ is a {\em
quasi-isomorphism} if $H^n\p$ is an isomorphism for all
$n\in\bbZ$. Note that two maps $\p,\psi\colon X\to Y$ induce the same
map $H^n\p=H^n\psi$, if $\p-\psi$ is null-homotopic.

The {\em derived category} $\bfD(\A)$ of $\A$ is obtained
from $\bfK(\A)$ by formally inverting all quasi-isomorphisms. To be
precise, one defines
$$\bfD(\A)=\bfK(\A)[S^{-1}]$$ as the localization of $\bfK(\A)$ with
respect to the class $S$ of all quasi-isomorphisms.  

\subsection{Localization}
Let $\C$ be a category and $S$ be a class of maps in $\C$. The {\em
localization} of $\C$ with respect to $S$ is a category $\C[S^{-1}]$,
together with a functor $Q\colon\C\to \C[S^{-1}]$ such that
\begin{enumerate}
\item[(L1)] $Q\s$ is an isomorphism for all $\s$ in $S$, and
\item[(L2)] any functor $F\colon\C\to\D$ such that $F\s$ is an
isomorphism for all $\s$ in $S$ factors uniquely through $Q$.
\end{enumerate}
Ignoring set-theoretic problems, one can show that such a localization
always exists. Let us give an explicit construction for
$\C[S^{-1}]$. To this end recall, that $S$ is a {\em multiplicative
system} if it satisfies the following conditions.
\begin{enumerate}
\item[(MS1)] If $\s,\t$ are composable maps in $S$, then $\t\comp \s$
is in $S$. The identity map $\id_X$ is in $S$ for all $X$ in $\C$.
\item[(MS2)] Let $\s\colon X\to Y$ be in $S$. Then every pair of maps
$Y'\to Y$ and $X\to X''$ in $\C$ can be completed to a pair of
commutative diagrams
$$\xymatrix{X'\ar[r]\ar[d]^{\s'}&X\ar[d]^{\s}&
X\ar[r]\ar[d]^{\s}&X''\ar[d]^{\s''}\\ Y'\ar[r]&Y&Y\ar[r]&Y''}$$ such
that $\s'$ and $\s''$ are in $S$.
\item[(MS3)] Let $\a,\b\colon X\to Y$ be maps in $\C$. Then there is
some $\s\colon X'\to X$ in $S$ with $\a\comp\s=\b\comp \s$ if and
only if there is some $\t\colon Y\to Y'$ in $S$ with $\t\comp\a=\t\comp\b$.
\end{enumerate}
Assuming that $S$ is a multiplicative system, one obtains the
following description of $\C[S^{-1}]$. The objects are those of
$\C$. Given objects $X$ and $Y$, the maps $X\to Y$ in $\C[S^{-1}]$ are
equivalence classes of diagrams
$$\xymatrix{X\ar[r]^\a&Y'&Y\ar[l]_\s}$$ with $\s$ in $S$, where two
pairs $(\a_1,\s_1)$ and $(\a_2,\s_2)$ are equivalent if there exists a
commutative diagram
$$\xymatrix{&Y_1\ar[d]\\ X\ar[r]^{\a_3}\ar[ru]^{\a_1}\ar[rd]_{\a_2}&Y_3&
Y\ar[lu]_{\s_1}\ar[ld]^{\s_2}\ar[l]_{\s_3}\\ &Y_2\ar[u]}$$ with $\s_3$ in
$S$. The composition of two pairs $(\a,\s)$ and $(\b,\t)$ is by
definition the pair $(\b'\comp\a,\s'\comp\t)$ where $\s'$ and $\b'$
are obtained from condition (MS2) as in the following commutative
diagram.
$$\xymatrix{&&Z''\\ &Y'\ar[ru]^{\b'}&& Z'\ar[lu]_{\s'}\\
X\ar[ru]^\a&&Y\ar[lu]_{\s}\ar[ru]^\b&&Z\ar[lu]_{\t}}$$ The universal
functor $\C\to\C[S^{-1}]$ sends a map $\a\colon X\to Y$ to the pair
$(\a,\id_Y)$. A pair $(\a,\s)$ is called a {\em fraction} because it is
identified  with $\s^{-1}\comp\a$ in $\C[S^{-1}]$.

\begin{exm} 
The quasi-isomorphisms in $\bfK(\A)$ form a multiplicative system; see
(\ref{ss:qis}).
\end{exm}

\subsection{An alternative definition}
Let $\A$ be an abelian category. 
We denote by $S$ the class of
quasi-isomorphisms in $\bfK(\A)$ and by $T$ the class of
quasi-isomorphisms in $\bfC(\A)$.
The derived category of $\A$ is by definition the localization of
$\bfK(\A)$ with respect $S$. Alternatively, one could take the
localization of $\bfC(\A)$ with respect $T$.

\begin{lem} 
The canonical functor $P\colon\bfC(\A)\to\bfK(A)$ induces a unique
isomorphism $\bar P$ making the following diagram of functors commutative.
$$\xymatrix{\bfC(\A)\ar[d]^P\ar[rr]^-{Q_T}&&\bfC(\A)[T^{-1}]\ar[d]^{\bar P}\\
\bfK(\A)\ar[rr]^-{Q_S}&&\bfK(\A)[S^{-1}]}$$
\end{lem}
\begin{proof} 
We have $T=P^{-1}(S)$. Thus $Q_S\comp P$ inverts all maps in $T$ and
induces therefore a unique functor $\bar P$ making the above diagram
commutative. On the other hand, $Q_T$ factors through $P$ via a
functor $F\colon \bfK(\A)\to\bfC(\A)[T^{-1}]$. This follows from the
fact that for two maps $\p,\psi\colon X\to Y$ in $\bfC(\A)$, we have
$Q_T\p=Q_T\psi$ if $\p-\psi$ is null-homotopic; see
\cite[Lemma~III.4.3]{GM}. The functor $F$ inverts all maps in $S$ and
factors therefore through $Q_S$ via a functor $\bar F\colon
\bfK(\A)[S^{-1}]\to \bfC(\A)[T^{-1}]$. Clearly, $Q_T\comp \bar F\comp
\bar P=Q_T$ and $Q_S\comp \bar P\comp\bar F=Q_S$. This implies $\bar F\comp
\bar P=\Id$ and $\bar P\comp\bar F=\Id$.
\end{proof}

\subsection{Extension groups}\label{ss:ext}
Let $\A$ be an abelian category. An object $A$ in $\A$ is identified
with the complex
$$\cdots\lto 0\lto A\lto 0\lto\cdots$$
concentrated in degree zero.
Given any complex $X$ in $\A$, we denote by $\Si X$ or $X[1]$ the
shifted complex with $$(\Si X)^n=X^{n+1}\quad\textrm{and}\quad
d^n_{\Si X}=-d^{n+1}_X.$$
Next we describe the maps in $\bfD(\A)$ for
certain complexes.
\begin{lem} 
  Let $Y$ be a complex in $\A$ such that each $Y^n$ is injective and
  $Y^n=0$ for $n\ll 0$.
\begin{enumerate}
\item Every quasi-isomorphism $\s\colon Y \to Y'$ has a left inverse
$\s'\colon Y'\to Y$ such that $\s'\comp\s=\id_Y$ in $\bfK(\A)$.
\item Given a complex $X$, the map
$\Hom_{\bfK(\A)}(X,Y)\to\Hom_{\bfD(\A)}(X,Y)$ is bijective.
\end{enumerate}
\end{lem}
\begin{proof} 
For (1), see (\ref{ss:Ktria}). To show (2), let
$X\stackrel{\a}\rightarrow Y'\stackrel{\s}\leftarrow Y$ be a diagram
representing a map $X\to Y$ in $\bfD(\A)$. Then $\s$ is a
quasi-isomorphism and has therefore a left inverse $\s'$ in $\bfK(\A)$, by
(1). Thus $(\a,\s)$ is equivalent to $(\s'\comp\a,\id_Y)$, which
belongs to the image of the canonical map
$\Hom_{\bfK(\A)}(X,Y)\to\Hom_{\bfD(\A)}(X,Y)$.

Now let $\a_1$ and $\a_2$ be maps $X\to Y$ such that $(\a_1,\id_Y)$
and $(\a_2,\id_Y)$ are equivalent. Then there is a quasi-isomorphism
$\s\colon Y\to Y'$ with $\s\comp\a_1=\s\comp\a_2$. The map $\s$ has a
left inverse in $\bfK(\A)$, by (1). Thus $\a_1=\a_2$.
\end{proof}

\begin{exm}
  Let $I$ be an injective object in $\A$. Then
$$\Hom_{\bfK(\A)}(-,I)\cong\Hom_\A(H^0-,I).$$
\end{exm}

A motivation for introducing the derived category is the fact that the
maps in $\bfD(\A)$ provide all homological information on objects in
$\A$.  This becomes clear from the following description of the
derived functors $\Ext_\A^n(-,-)$.  Here, we use the convention that
$\Ext_\A^n(-,-)$ vanishes for $n<0$.

\begin{lem} 
For all objects $A,B$ in $\A$ and $n\in\bbZ$, there is a canonical
isomorphism
$$\Ext^n_\A(A,B)\lto\Hom_{\bfD(\A)}(A,\Si^nB).$$
\end{lem}
\begin{proof}
We give a proof in case $\A$ has enough injective objects; otherwise
see \cite[III.3.2]{V}. Choose for $B$ an injective resolution $\bfi
B$. Thus we have a quasi-isomorphism $B\to \bfi B$ which induces the
following isomorphism.
\begin{align*}
\Ext^n_\A(A,B)&= H^n\Hom_{\A}(A,\bfi B)
\cong\Hom_{\bfK(\A)}(A,\Si^n\bfi B)\\
&\cong \Hom_{\bfD(\A)}(A,\Si^n\bfi B)
\cong\Hom_{\bfD(\A)}(A,\Si^nB)
\end{align*}
\end{proof}
Note as a consequence that the canonical functor $\A\to\bfD(\A)$ which
sends an object to the corresponding complex concentrated in degree
zero is fully faithful. In fact, it identifies $\A$ with the full
subcategory of complexes $X$ in $\bfD(\A)$ such that $H^nX=0$ for all
$n\neq 0$.

\subsection{Hereditary categories} 
Let $\A$ be a {\em hereditary} abelian category, that is,
$\Ext_\A^2(-,-)$ vanishes. In this case, there is an explicit
description of all objects and maps in $\bfD(\A)$. Every complex $X$
is completely determined by its cohomology because there is an
isomorphism between $X$ and 
$$\cdots \lto H^{n-1}X\lto[0]H^nX\lto[0]H^{n+1}X\lto\cdots.$$ To
construct this isomorphism, note that the vanishing of
$\Ext^2_\A(H^nX,-)$ implies the existence of a commutative diagram
$$\xymatrix{0\ar[r]&X^{n-1}\ar[r]\ar[d]&E^n\ar[r]\ar[d]&
H^nX\ar[r]\ar@{=}[d]&0\\ 0\ar[r]&\Im d^{n-1}\ar[r]&\Ker
d^n\ar[r]&H^nX\ar[r]&0}$$ with exact rows. We obtain the following
commutative diagram
$$\xymatrix{ \cdots\ar[r]&0\ar[r]&0\ar[r]&H^nX\ar[r]&0\ar[r]&\cdots\\
\cdots\ar[r]&0\ar[u]\ar[r]\ar[d]&X^{n-1}\ar[u]\ar[r]\ar@{=}[d]&E^n\ar[u]\ar[r]\ar[d]&
0\ar[u]\ar[r]\ar[d]&\cdots\\
\cdots\ar[r]&X^{n-2}\ar[r]&X^{n-1}\ar[r]&X^{n}\ar[r]&X^{n+1}\ar[r]&\cdots
}$$ and the vertical maps induce cohomology isomorphism in degree $n$.
Thus we have in $\bfD(\A)$
$$\coprod_{n\in\bbZ}\Si^{-n}(H^nX)\cong X\cong
\prod_{n\in\bbZ}\Si^{-n}(H^nX).$$ Using the description of
$\Hom_{\bfD(\A)}(A,B)$ for objects $A,B$ in $\A$, we see that each map
$X\to Y$ in $\bfD(\A)$ is given by a family of maps in
$\Hom_\A(H^nX,H^nY)$ and a family of extensions in $\Ext_\A^1(H^nX,
H^{n-1}Y)$, with $n\in\bbZ$. Thus $$\bfD(\A)=\bigsqcup_{n\in\bbZ}\Si^n\A$$ 
with non-zero maps $\Si^i\A\to\Si^j\A$ only if $j-i\in\{0,1\}$.

\begin{exm} 
The category  of abelian groups is hereditary. More
generally, a module category is hereditary if and only if every
submodule of a projective module is projective.
\end{exm}

\subsection{Bounded derived categories}
Let $\A$ be an additive category. Consider the following full
subcategories of $\bfC(\A)$.
\begin{align*}
\bfC^-(\A)&=\{X\in\bfC(\A)\mid X^n=0\textrm{ for }n\gg 0\}\\
\bfC^+(\A)&=\{X\in\bfC(\A)\mid X^n=0\textrm{ for }n\ll 0\}\\
\bfC^b(\A)&=\{X\in\bfC(\A)\mid X^n=0\textrm{ for }|n|\gg 0\}\\
\end{align*}
For $*\in\{-,+,b\}$, let the homotopy category $\bfK^*(\A)$ be
the quotient of $\bfC^*(\A)$ modulo null-homotopic maps, and let the
derived category $\bfD^*(\A)$ be the localization with respect to all
quasi-isomorphisms.

\begin{lem} 
For each $*\in\{-,+,b\}$, the inclusion $\bfC^*(\A)\to\bfC(\A)$ induces
fully faithful functors $\bfK^*(\A)\to\bfK(\A)$ and $\bfD^*(\A)\to\bfD(\A)$.
\end{lem}
\begin{proof} See \cite[III.1.2.3]{V}.
\end{proof}

\subsection{Notes} 
Derived categories were introduced by Grothendieck and his school in
the 1960s; see Verdier's (posthumously published) th\`ese \cite{V}.
The calculus of fractions which appears in the description of the
localization of a category was developed by Gabriel and Zisman
\cite{GZ}. For a modern treatment of derived categories and related
material, see \cite{GM,KS,W}.

\section{Triangulated categories}

The derived category $\bfD(\A)$ of an abelian category is an additive
category. There is some additional structure which complements the
abelian structure of $\bfC(\A)$. The axiomatization of this structure leads
to the notion of a triangulated category.

\subsection{The axioms}\label{ss:axioms}

Let $\T$ be an additive category with an equivalence
$\Si\colon\T\to\T$. A {\em triangle} in $\T$ is a sequence
$(\a,\b,\g)$ of maps
$$X\lto[\a] Y\lto[\b] Z\lto[\g]\Si X,$$ and a morphism between two
triangles $(\a,\b,\g)$ and $(\a',\b',\g')$ is a triple $(\p_1,\p_2,\p_3)$
of maps in $\T$ making the following diagram commutative.
$$\xymatrix{
X\ar[r]^\a\ar[d]^{\p_1}&Y\ar[r]^\b\ar[d]^{\p_2}&Z\ar[r]^\g\ar[d]^{\p_3}&\Si
X\ar[d]^{\Si\p_1}\\ X'\ar[r]^{\a'}&Y'\ar[r]^{\b'}&Z'\ar[r]^{\g'}&\Si X'
}$$ The category $\T$ is called {\em triangulated} if it is equipped
with a class of distinguished triangles (called {\em exact triangles})
satisfying the following conditions.
\begin{enumerate}
\item[(TR1)] A triangle isomorphic to an exact triangle is exact. For
each object $X$, the triangle $0\to X\xto{\id} X\to 0$ is exact.
Each map $\a$ fits into an exact triangle $(\a,\b,\g)$.
\item[(TR2)] A triangle $(\a,\b,\g)$ is exact if and only if
$(\b,\g,-\Si\a)$ is exact.
\item[(TR3)] Given two exact triangles $(\a,\b,\g)$ and
$(\a',\b',\g')$, each pair of maps $\p_1$ and $\p_2$ satisfying
$\p_2\comp\a=\a'\comp\p_1$ can be completed to a morphism
$$\xymatrix{
X\ar[r]^\a\ar[d]^{\p_1}&Y\ar[r]^\b\ar[d]^{\p_2}&Z\ar[r]^\g\ar[d]^{\p_3}&\Si
X\ar[d]^{\Si\p_1}\\ X'\ar[r]^{\a'}&Y'\ar[r]^{\b'}&Z'\ar[r]^{\g'}&\Si X'
}$$ of triangles.
\item[(TR4)] Given exact triangles $(\a_1,\a_2,\a_3)$,
$(\b_1,\b_2,\b_3)$, and $(\g_1,\g_2,\g_3)$ with $\g_1=\b_1\comp\a_1$,
there exists an exact triangle $(\d_1,\d_2,\d_3)$ making
the following diagram commutative.
$$\xymatrix{X\ar[r]^{\a_1}\ar@{=}[d]&Y\ar[r]^{\a_2}\ar[d]^{\b_1}&
U\ar[r]^{\a_3}\ar[d]^{\d_1}& \Si X\ar@{=}[d]\\
X\ar[r]^{\g_1}&Z\ar[r]^{\g_2}\ar[d]^{\b_2}&
V\ar[r]^{\g_3}\ar[d]^{\d_2}&\Si X\ar[d]^{\Si\a_1}\\
&W\ar@{=}[r]\ar[d]^{\b_3}& W\ar[d]^{\d_3}\ar[r]^{\b_3}&\Si Y\\ &\Si
Y\ar[r]^{\Si\a_2}&\Si U }$$
\end{enumerate}

The category $\T$ is called {\em pre-triangulated} if the axioms (TR1)
-- (TR3) are satisfied.

\subsection{The octahedral axiom}
Let $\T$ be a pre-triangulated category.  The axiom (TR4) is known as
{\em octahedral axiom} because the four exact triangles can be
arranged in a diagram having the shape of an octahedron. The exact
triangles $A\to B\to C\to\Si A$ are represented by faces of the form
$$\xymatrix{&C\ar[ld]|-{+}\\A\ar[rr]&&B\ar[lu]}$$
and the other four
faces are commutative triangles.

Let us give a more intuitive formulation of the octahedral axiom which
is based on the notion of a homotopy cartesian square.  Call a
commutative square
$$\xymatrix{X\ar[r]^{\a'}\ar[d]^{\a''}&Y'\ar[d]^{\b'}\\Y''
\ar[r]^{\b''}&Z}$$ {\em homotopy cartesian} if there exists an exact
triangle
$$X\xto{\smatrix{\a'\\ \a''}}Y'\amalg
Y''\xto{\smatrix{\b'&-\b''}}Z\lto[\g]\Si X.$$ The map $\g$ is called a
{\em differential} of the homotopy cartesian square.  Note that a
differential of the homotopy cartesian square changes its sign if the
square is flipped along the main diagonal.
\begin{enumerate}
\item[(TR4')] Every pair of maps $X\to Y$ and $X\to X'$ can be
completed to a morphism
$$\xymatrix{
  X\ar[r]\ar[d]^{}&Y\ar[r]\ar[d]^{}&Z\ar[r]\ar@{=}[d]^{}&\Si
  X\ar[d]^{}\\ X'\ar[r]^{}&Y'\ar[r]^{}&Z\ar[r]^{}&\Si X' }$$ between
  exact triangles such that the left hand square is homotopy cartesian
  and the composite $Y'\to Z\to\Si X$ is a differential.
\end{enumerate}

We can think of a homotopy cartesian square as the triangulated
analogue of a pull-back and push-out square in an abelian
category. Recall that a square
$$\xymatrix{X\ar[r]^{\a'}\ar[d]^{\a''}&Y'\ar[d]^{\b'}\\Y''
  \ar[r]^{\b''}&Z}$$
is a pull-back and a push-out if and only if the
sequence
$$0\lto X\xto{\smatrix{\a'\\ \a''}}Y'\amalg
Y''\xto{\smatrix{\b'&-\b''}}Z\lto 0$$ is exact. The axiom (TR4') is
the triangulated analogue of the fact that parallel maps in a
pull-back square have isomorphic kernels, whereas parallel maps in a
push-out square have isomorphic cokernels.

\begin{prop} 
  The axioms (TR4) and (TR4') are equivalent for any pre-triangulated category.
\end{prop}

The proof can be found in the appendix.

\subsection{Cohomological functors}\label{ss:hom}
Let $\T$ be a pre-triangulated category.  Note that the axioms of a
pre-triangulated category are symmetric in the sense that the opposite
category $\T^\op$ carries a canonical pre-triangulated structure.

Given an abelian category $\A$, a functor $\T\to\A$ is called {\em
cohomological} if it sends each exact triangle in $\T$ to an exact
sequence in $\A$.

\begin{lem}\label{le:cohfun}
 For each object $X$ in $\T$, the representable functors
$$\Hom_\T(X,-)\colon\T\lto\Ab\quad\textrm{and}\quad
\Hom_\T(-,X)\colon\T^\op\lto\Ab$$
into the category $\Ab$ of abelian groups are cohomological functors.
\end{lem}
\begin{proof} 
  Fix an exact triangle $U\to[\a]V\to[\b]W\to[\g]\Si U$. We need to
  show the exactness of the induced sequence
$$\Hom_\T(X,U)\lto\Hom_\T(X,V)\lto\Hom_\T(X,W)\lto\Hom_\T(X,\Si U).$$
It is sufficient to check exactness at one place, say $\Hom_\T(X,V)$,
by (TR2).  To this end fix a map $\p\colon X\to V$ and consider the
following diagram.
$$\xymatrix{X\ar[r]^\id&X\ar[r]\ar[d]^\p&0\ar[r]&\Si X\\
U\ar[r]^\a&V\ar[r]^\b&W\ar[r]^\g&\Si U}$$ If $\p$ factors through
$\a$, then (TR3) implies the existence of a map $0\to W$ making the
diagram commutative. Thus $\b\comp\p=0$. Now assume
$\b\comp\p=0$. Applying (TR2) and (TR3), we find a map $X\to U$ making
the diagram commutative. Thus $\p$ factors through $\a$. 
\end{proof}

\subsection{Uniqueness of exact triangles}\label{ss:uni}

Let $\T$ be a pre-triangulated category. Given a map $\a\colon X\to Y$
in $\T$ and two exact triangles $\Delta=(\a,\b,\g)$ and
$\Delta'=(\a,\b',\g')$ which complete $\a$, there exists a comparison
map $(\id_X,\id_Y,\p)$ between $\Delta$ and $\Delta'$, by (TR3). The
map $\p$ is an isomorphism, by the following lemma, but it need not to
be unique.

\begin{lem}
Let $(\p_1,\p_2,\p_3)$ be a morphism between exact triangles.  If two
maps from $\{\p_1,\p_2,\p_3\}$ are isomorphisms, then also the third.
\end{lem}
\begin{proof} Use lemma~(\ref{ss:hom}) and apply the 5-lemma.
\end{proof}

The third object $Z$ in an exact triangle $X\to Y\to Z\to\Si X$ is
called the {\em cone} or the {\em cofiber} of the map $X\to Y$. 

\subsection{$\bfK(\A)$ is triangulated}\label{ss:Ktria}

Let $\A$ be an additive category and let $\bfK(\A)$ be the homotopy
category of complexes. Consider the equivalence
$\Si\colon\bfK(\A)\to\bfK(\A)$ which takes a complex to its shifted
complex. Given a map $\a\colon X\to Y$ of complexes, the {\em mapping
cone} is the complex $Z$ with $Z^n=X^{n+1}\amalg Y^n$ and differential
$\smatrix{-d^{n+1}_X&0\\ \a^{n+1}&d^n_Y}$. The mapping cone fits into
a {\em mapping cone sequence} $$X\lto[\a]Y\lto[\b]Z\lto[\g]\Si X$$
which is defined in degree $n$ by the following sequence.
$$X^n\lto[\a^n] Y^n\lto[\smatrix{0\\ \id}] X^{n+1}\amalg Y^{n}
\lto[\smatrix{-\id&0}] X^{n+1}$$ By definition, a triangle in
$\bfK(\A)$ is {\em exact} if it is isomorphic to a mapping cone
sequence as above. It is easy to verify the axioms (TR1) -- (TR4); see
\cite[II.1.3.2]{V} or (\ref{ss:stab}).  Thus $\bfK(\A)$ is a
triangulated category.

Any mapping cone sequence of complexes induces a long exact sequence
when one passes to its cohomology. To make this precise, we identify
$H^i(\Si^jX)=H^{i+j}X$ for every complex $X$ and all $i,j$.
\begin{lem}\label{le:longex}
Let $\A$ be an abelian category. An exact triangle
$$X\lto[\a]Y\lto[\b]Z\lto[\g]\Si X$$ in $\bfK(\A)$ induces the following long
exact sequence.
$$ \cdots\lto
H^{n-1}Z\lto[H^{n-1}\g]H^{n}X\lto[H^{n}\a]H^{n}Y\lto[H^{n}\b]H^{n}Z\lto
[H^{n}\g]H^{n+1}X\lto\cdots$$
\end{lem}
\begin{proof}
We may assume that the triangle is a mapping cone sequence as
above. The short exact sequence $0\to Y\to Z\to\Si X\to 0$ of
complexes induces a long exact sequence
$$ \cdots\lto H^{n-1}(\Si
X)\lto[\d^{n-1}]H^{n}Y\lto[H^{n}\b]H^{n}Z\lto[H^{n}\g]H^{n}(\Si X)\lto
[\d^{n}]H^{n+1}Y\lto\cdots$$ with connecting morphism $\d^*$. This
follows from the Snake Lemma. Now observe that $\d^n=H^{n+1}\a$.
\end{proof}

\begin{rem}
  The triangulated structure of $\bfK(\A)$ allows a quick proof of the
  following fact. Given a complex $Y$ such that each $Y^n$ is
  injective and $Y^n=0$ for $n\ll 0$, every quasi-isomorphism
  $\b\colon Y\to Z$ has a left inverse $\b'$ such that $\b'\comp\b=\id_Y$
  in $\bfK(\A)$. In order to construct $\b'$, complete $\b$ to an
  exact triangle $X\to[\a] Y\to[\b] Z\to[\g]\Si X$. Then $\a$ is
  null-homotopic since $H^nX=0$ for all $n\in\bbZ$. Thus $\id_Y$
  factors through $\b$, by lemma~(\ref{ss:hom}).
\end{rem}

\subsection{Notes}

Triangulated categories were introduced independently in algebraic
geometry by Verdier in his th\`ese \cite{V}, and in algebraic topology
by Puppe \cite{P}. The basic properties of triangulated categories can
be found in \cite{V}.  The reformulation (TR4') of the octahedral
axiom (TR4) which is given here is a variation of a reformulation due
to Dlab, Parshall, and Scott \cite{PS}. There is an equivalent
formulation of (TR4) due to May which displays the exact triangles in
a diagram having the shape of a braid \cite{M}.  Another reformulation
is discussed in Neeman's book \cite{N3}, which also contains a
discussion of homotopy cartesian squares. We refer to \cite{KS,N3} for a
modern treatment of triangulated categories.

\section{Localization of triangulated categories}

The triangulated structure of $\bfK(\A)$ induces a triangulated
structure of $\bfD(\A)$ via the localization functor
$\bfK(\A)\to\bfD(\A)$. This follows from the fact that the
quasi-isomorphisms form a class of maps in $\bfK(\A)$ which is
compatible with the triangulation.

\subsection{Quasi-isomorphisms}\label{ss:qis}

Let $\T$ be a triangulated category and $S$ be a class of maps which
is a multiplicative system. Then $S$ is {\em compatible with the
triangulation} if
\begin{enumerate}
\item[(MS4)] given $\s$ in $S$, the map $\Si^n\s$ belongs to $S$ for
all $n\in\bbZ$, and
\item[(MS5)] given a map $(\p_1,\p_2,\p_3)$ 
between exact triangles with $\p_1$ and $\p_2$ in $S$, there is also a
map $(\p_1,\p_2,\p'_3)$ with $\p'_3$ in $S$.
\end{enumerate}

\begin{lem}\label{le:comp}
Let $H\colon\T\to\A$ be a cohomological functor. Then the class $S$ of
maps $\s$ in $\T$ such that $H(\Si^n\s)$ is an isomorphism for all
$n\in\bbZ$ form a multiplicative system, compatible with the
triangulation.
\end{lem}
\begin{proof}
(MS1) and (MS4) are immediate consequences of the definition. (MS5)
follows from the 5-lemma. To show (MS2), let $\s\colon X\to Y$ be in
$S$ and $Y'\to Y$ be an arbitrary map. Complete $Y'\to Y$ to an exact
triangle. Applying (TR2) and (TR3), we obtain the following morphism
between exact triangles
$$\xymatrix{
X'\ar[r]\ar[d]^{\s'}&X\ar[r]\ar[d]^{\s}&Y''\ar[r]\ar@{=}[d]&\Si
X'\ar[d]^{\Si\s'}\\ Y'\ar[r]&Y\ar[r]&Y''\ar[r]&\Si Y' }$$ and the
5-lemma shows that $\s'$ belongs to $S$. It remains to check (MS3).
Let $\a,\b\colon X\to Y$ be maps in $\T$ and $\s\colon X'\to X$ in $S$
such that $\a\comp\s=\b\comp \s$. Complete $\s$ to an exact triangle
$X'\to[\s] X\to[\p] X''\to\Si X'$. Then $\a-\b$ factors through $\p$
via some map $\psi\colon X''\to Y$. Now complete $\psi$ to an exact
triangle $X''\to[\psi] Y\to[\t] Y'\to\Si X''$. Then 
$\t$ belongs to $S$ and $\t\comp\a=\t\comp\b$.
\end{proof}
\subsection{$\bfD(\A)$ is triangulated}

An {\em exact} functor $\T\to\U$ between triangulated categories is a
pair $(F,\eta)$ consisting of a functor $F\colon\T\to\U$ and a natural
isomorphism $\eta\colon F\comp \Si_\T\to\Si_\U\comp F$ such that for
every exact triangle $X\to[\a]Y\to[\b]Z\to[\g]\Si X$ in $\T$ the
triangle
$$FX\lto[F\a]FY\lto[F\b]FZ\xto{\eta_X\comp F\g}\Si (FX)$$ is exact in $\U$.

\begin{lem}\label{le:qis}
Let $\T$ be a triangulated category and $S$ be a multiplicative system
of maps which is compatible with the triangulation. Then the
localization $\T[S^{-1}]$ carries a unique triangulated structure such
that the canonical functor $\T\to\T[S^{-1}]$  is exact.
\end{lem}
\begin{proof}
The equivalence $\Si\colon\T\to\T$ induces a unique equivalence
$\T[S^{-1}]\to\T[S^{-1}]$ which commutes with the canonical functor
$Q\colon\T\to\T[S^{-1}]$ . This follows from (MS4). Now take as exact
triangles in $\T[S^{-1}]$ all those isomorphic to images of exact
triangles in $\T$. It is straightforward to verify the axioms (TR1) --
(TR4); see \cite[II.2.2.6]{V}. The functor $Q$ is exact by
construction.
\end{proof}

Let $\A$ be an abelian category. The quasi-isomorphisms in $\bfK(\A)$
form a multiplicative system which is compatible with the
triangulation, by lemma~(\ref{le:comp}). In fact,
$H^0\colon\bfK(\A)\to\A$ is a cohomological functor by
lemma~(\ref{le:longex}), and a map $\s$ in $\bfK(\A)$ is a
quasi-isomorphism if and only if $H^0(\Si^n\s)=H^n\s$ is an
isomorphism for all $n\in\bbZ$.  Thus $\bfD(\A)$ is triangulated and
the canonical functor $\bfK(\A)\to\bfD(\A)$ is exact.

The canonical functor $\A\to \bfD(\A)$ sends every exact sequence 
$$E\colon 0\lto A\lto[\a] B\lto[\b] C\lto 0$$
in $\A$ to an exact
triangle
$$A\lto[\a] B\lto[\b] C\lto[\g] \Si A$$ in $\bfD(\A)$, where $\g$
denotes the map corresponding to $E$ under the canonical
isomorphism
$$\Ext_\A^1(C,A)\cong\Hom_{\bfD(\A)}(C,\Si A).$$ In fact, the mapping
cone of $\a$ produces the following exact triangle
$$\xymatrix{
\cdots\ar[r]&0\ar[r]\ar[d]&0\ar[r]\ar[d]&A\ar[r]\ar[d]^\a&
0\ar[r]\ar[d]&\cdots\\
\cdots\ar[r]&0\ar[r]\ar[d]&0\ar[r]\ar[d]&B\ar[r]\ar[d]^\id&
0\ar[r]\ar[d]&\cdots\\
\cdots\ar[r]&0\ar[r]\ar[d]&A\ar[r]^\a\ar[d]^{-\id}&B\ar[r]\ar[d]&
0\ar[r]\ar[d]&\cdots\\
\cdots\ar[r]&0\ar[r]&A\ar[r]&0\ar[r]&0\ar[r]&\cdots}$$ and the map
$\b$ induces a quasi-isomorphism between the mapping cone of $\a$ and
the complex corresponding to $C$.

\begin{exm}
  Let $F\colon\A\to\B$ be an additive functor between additive
  categories. Then $F$ induces an exact functor $\bfK(F)\colon
  \bfK(\A)\to\bfK(\B)$.  Now suppose that $\A$ and $\B$ are abelian
  and $F$ is exact.  Then $\bfK(F)$ sends quasi-isomorphisms to
  quasi-isomorphisms and induces therefore an exact functor
  $\bfD(\A)\to\bfD(\B)$.
\end{exm}

\subsection{Triangulated and thick subcategories}\label{ss:tria}
Let $\T$ be a triangulated category. A non-empty full subcategory $\S$ is a
{\em triangulated subcategory} if the following conditions hold.
\begin{enumerate}
\item[(TS1)] $\Si^n X\in\S$ for all $X\in\S$ and $n\in\bbZ$.
\item[(TS2)] Let $X\to Y\to Z\to\Si X$ be an exact triangle in $\T$.
  If two objects from $\{X,Y,Z\}$ belong to $\S$, then also the third.
\end{enumerate}
A triangulated subcategory $\S$ is {\em thick} if in addition the following
condition holds.
\begin{enumerate}
\item[(TS3)] Every direct factor of an object in $\S$ belongs to $\S$,
that is, a decomposition $X=X'\amalg X''$ for $X\in\S$ implies
$X'\in\S$.
\end{enumerate}
Note that a triangulated subcategory $\S$ inherits a canonical
triangulated structure from $\T$. 

Given a class $\S_0$ of objects in $\T$, one can construct inductively
the triangulated subcategory {\em generated by} $\S_0$ as
follows. Denote for two classes $\U$ and $\V$ of objects in $\T$ by
$\U * \V$ the class of objects $X$ occuring in an exact triangle $U\to
X\to V\to \Si U$ with $U\in\U$ and $V\in\V$. Note that the operation
$*$ is associative by the octahedral axiom. Now let $\S_1$ be the
class of all $\Si^nX$ with $X\in\S_0$ and $n\in\bbZ$.  For $r>0$, let
$\S_r=\S_1*\S_1*\cdots *\S_1$ be the product with $r$ factors.

\begin{lem}
Let $\S_0$ be a class of objects in $\T$.
\begin{enumerate} 
\item  The full subcategory of objects in
$\S=\bigcup_{r\geqslant 0}\S_r$ is the smallest full triangulated
subcategory of $\T$ which contains $\S_0$.
\item The full subcategory of direct factors of objects in $\S$ is the
smallest full thick subcategory of $\T$ which contains $\S_0$.
\end{enumerate}
\end{lem}
\begin{proof}
  (1) is clear. To prove (2), we need to show that all direct factors
  of objects in $\S$ form a triangulated subcategory. To this end fix
  an exact triangle
  $$\Delta'\colon X'\lto[\a'] Y'\lto[\b'] Z'\lto[\g']\Si X'$$ such
  that $X=X'\amalg X''$ and $Y=Y'\amalg Y''$ belong to $\S$. We need
  to show that $Z'$ is a direct factor of some object in $\S$.
  Complete $\a=\smatrix{\a'&0\\0&0}$ to an exact triangle
  $$\Delta\colon X\lto[\a] Y\lto[\b] Z\lto[\g]\Si X$$ with $Z$ in
  $\S$. We have obvious maps $\iota\colon\Delta'\to\Delta$ and
  $\pi\colon\Delta\to\Delta'$ such that $\pi\comp\iota$ is invertible,
  by lemma~(\ref{ss:uni}).  Thus $Z'$ is a direct factor of $Z$.
\end{proof}

\begin{exm}
(1) A complex $X$ in some abelian category $\A$ is {\em acyclic} if
$H^nX=0$ for all $n\in\bbZ$. The acyclic complexes form a thick
subcategory in $\bfK(\A)$. A map between complexes is a
quasi-isomorphism if and only if its mapping cone is acyclic.  The
canonical functor $\bfK(\A)\to\bfD(\A)$ annihilates a map in
$\bfK(\A)$ if and only if it factors through an acyclic complex; see
(\ref{se:ker}).

(2)  The bounded derived category $\bfD^b(\A)$ of an abelian category
  $\A$ is generated as a triangulated category by the objects in $\A$,
  viewed as complexes concentrated in degree zero.

(3) Let $\A$ be an abelian category with enough injective objects.
Then $\bfD^b(\A)$ is generated as a triangulated category by all
injective objects, if and only if every object in $\A$ has finite
injective dimension.
\end{exm}

\subsection{The kernel of a localization}\label{se:ker}
Let $\T$ be a triangulated category and let $F\colon\T\to\U$ be an
additive functor. The {\em kernel} $\Ker F$ of $F$ is by definition
the full subcategory of $\T$ which is formed by all objects $X$ such
that $FX=0$. If $F$ is an exact functor into a triangulated category,
then $\Ker F$ is a thick subcategory of $\T$. Also, if $F$ is a
cohomological functor into an abelian category, then
$\bigcap_{n\in\bbZ}\Ker (F\comp\Si^n)$ is a thick subcategory of $\T$.

Next we assume that $F$ is a localization functor and describe the
maps $\a$ in $\T$ such that $F\a=0$. Note that $FX=0$ for an object
$X$ if and only if $F\id_X=0$.
\begin{lem}\label{le:ker}
Let $\T$ be a triangulated category and $S$ be a multiplicative system
of maps in $\T$ which is compatible with the triangulation. The
following are equivalent for a map $\a\colon X\to Y$ in $\T$.
\begin{enumerate}
\item The canonical functor $Q\colon\T\to\T[S^{-1}]$ annihilates $\a$.
\item The exists a map $\s\colon Y\to Z$ in $S$ with $\s\comp\a=0$.
\item The map $\a$ factors through the cone of a map in $S$.
\end{enumerate}
\end{lem}
\begin{proof} 
(1) $\Leftrightarrow$ (2): The functor $Q$ sends $\a$ to the fraction
$(\a,\id_Y)$. Thus $Q$ annihilates $\a$ if and only if the fractions
$(\a,\id_Y)$ and $(0,\id_Y)$ are equivalent.  Now observe that both
fractions are equivalent if and only if $\s\comp\a=0$ for some $\s\in
S$.

(2) $\Leftrightarrow$ (1): Use
lemma~(\ref{le:cohfun}).
\end{proof}

\subsection{Verdier localization}
Let $\T$ be a triangulated category. Given a triangulated subcategory
$\S$, we denote by $S(\S)$ the class of maps in $\T$ such that the
cone of $\a$ belongs to $\S$.
 
\begin{lem}
Let $\T$ be a triangulated category and $\S$ be a triangulated
subcategory. Then $S(\S)$ is a multiplicative system
which is compatible with the triangulation of $\T$.
\end{lem}
\begin{proof} See \cite[II.2.1.8]{V}
\end{proof}

The {\em Verdier localization} of $\T$ with respect to a triangulated
subcategory $\S$ is by definition the localization
$$\T/\S=\T[S(\S)^{-1}]$$ together with the canonical functor
$\T\to\T/\S$.

\begin{prop}
Let $\T$ be a triangulated category and $\S$ a full triangulated
subcategory. Then the category $\T/\S$ and the canonical functor
$Q\colon\T\to\T/\S$ have the following properties.
\begin{enumerate}
\item The category $\T/\S$ carries a unique triangulated structure
such that $Q$ is exact.
\item The kernel $\Ker Q$ is the smallest thick subcategory containing
$\S$.
\item Every exact functor $\T\to\U$ annihilating $\S$ factors
uniquely through $Q$ via an exact functor $\T/\S\to\U$.
\item Every cohomological functor $\T\to\A$ annihilating $\S$ factors
uniquely through $Q$ via a cohomological functor $\T/\S\to\A$.
\end{enumerate}
\end{prop}
\begin{proof}
(1) follows from lemma~(\ref{le:qis}) and (2) from lemma~(\ref{le:ker}).

(3) An exact functor $F\colon\T\to\U$ annihilating $\S$ inverts every
    map in $S(\S)$. Thus there exists a unique functor $\bar
    F\colon\T/\S\to\U$ such that $F=\bar F\comp Q$. The functor $\bar
    F$ is exact because an exact triangle $\Delta$ in $\T/\S$ is up to
    isomorphism of the form $Q\Gamma$ for some exact triangle $\Gamma$
    in $\T$. Thus $\bar F\Delta=F\Gamma$ is exact.

(4) Analogous to (3).
\end{proof}

\subsection{Notes}

Localizations of triangulated categories are discussed in Verdier's
th\`ese \cite{V}. In particular, he introduced the localization
$\T/\S$ of a triangulated category $\T$ with respect to a triangulated
subcategory $\S$.

\section{Brown representability}

The Brown representability theorem provides a useful characterization
of the representable functors $\Hom_\T(-,X)$ for a class of
triangulated categories $\T$. The proof is based on a universal
embedding of $\T$ into an abelian category which is of independent
interest.  Note that such an embedding reverses the direction of the
construction of the derived category which provides an embedding of an
abelian category into a triangulated category.

\subsection{Coherent functors}\label{ss:coh}

Let $\A$ be an additive category.  We consider functors
$F\colon\A^\op\to\Ab$ into the category of abelian groups and call a
sequence $F'\to F\to F''$ of functors {\em exact} if the induced
sequence $F'X\to FX\to F''X$ of abelian groups is exact for all $X$ in
$\A$. A functor $F$ is said to be {\em coherent} if there exists an
exact sequence (called {\em presentation})
$$\A(-,X)\longrightarrow \A(-,Y)\longrightarrow F\longrightarrow 0.$$
Here, we simplify our notation and write $\A(X,Y)$ for the set of maps
$X\to Y$.  The natural transformations between two coherent functors
form a set by Yoneda's lemma, and the coherent functors $\A^\op\to\Ab$
form an additive category with cokernels. We denote this category by
$\widehat\A$. A basic tool is the fully faithful {\em Yoneda functor}
$$\A\lto\widehat\A,\quad X\mapsto \A(-,X).$$
Recall that a map $X\to Y$ is a {\em weak kernel} for $Y\to Z$
if the induced sequence
$$\A(-,X)\lto \A(-,Y)\lto \A(-,Z)$$ is exact.

\begin{lem} 
\begin{enumerate}
\item If $\A$ has weak kernels, then $\widehat \A$ is an abelian
category.
\item If $\A$ has arbitrary coproducts, then $\widehat \A$ has
arbitrary coproducts and the Yoneda functor preserves all coproducts.
\end{enumerate}
\end{lem}
\begin{proof} 
(1) The category $\widehat\A$ has cokernels, and it is therefore
sufficient to show that $\widehat\A$ has kernels. To this end fix a
map $F_1\to F_2$ with the following presentation.
$$\xymatrix{\A(-,X_1)\ar[r]\ar[d]&\A(-,Y_1)\ar[r]\ar[d]&F_1\ar[r]\ar[d]&0\\
\A(-,X_2)\ar[r]&\A(-,Y_2)\ar[r]&F_2\ar[r]&0}$$ We construct the kernel
$F_0\to F_1$ by specifying the following presentation.
$$\xymatrix{\A(-,X_0)\ar[r]\ar[d]&\A(-,Y_0)\ar[r]\ar[d]&F_0\ar[r]\ar[d]&0\\
\A(-,X_1)\ar[r]&\A(-,Y_1)\ar[r]&F_1\ar[r]&0}$$
First the map $Y_0\to Y_1$ is obtained from the weak
kernel sequence
$$Y_0\lto X_2\amalg Y_1\lto Y_2.$$
Then the maps $X_0\to
X_1$ and $X_0\to Y_0$ are obtained from the weak kernel sequence
$$X_0\lto X_1\amalg Y_0\lto Y_1.$$

(2) For every family of functors $F_i$ having a presentation
$$\A(-,X_i)\lto[(-,\p_i)] \A(-,Y_i)\lto F_i\lto 0,$$
the coproduct $F=\coprod_i F_i$ has a presentation
$$\A(-,\coprod_iX_i)\lto[(-,\amalg\p_i)] \A(-,\coprod_iY_i)\lto F\lto 0.$$
Thus coproducts in $\widehat\A$ are not computed pointwise.
\end{proof}

\begin{exm}
(1) Let $\A$ be an abelian category and suppose $\A$ has enough
    projective objects. Let $\P$ denote the full subcategory of
    projective objects in $\A$. Then the functor
$$\A\lto\widehat\P,\quad X\mapsto\Hom_\A(-,X)|_\P,$$ is an equivalence.

(2)  Let $X\to Y\to Z\to\Si X$ be an exact triangle in some triangulated
  category. Then $X\to Y$ is a weak kernel of $Y\to Z$; see
  lemma~(\ref{ss:hom}).
\end{exm}

\subsection{The abelianization of a triangulated category}\label{ss:abel}

Let $\T$ be a triangulated category. The Yoneda functor
$\T\to\widehat\T$ is the universal cohomological functor for $\T$.

\begin{lem}
\begin{enumerate}
\item The category $\widehat\T$ is abelian and the Yoneda functor
$H_\T\colon \T\to\widehat\T$ is cohomological.
\item Let $\A$ be an abelian category and $H\colon\T\to\A$ be a
cohomological functor. Then there is (up to a unique isomorphism) a
unique exact functor $\bar H\colon\widehat\T\to\A$ such that $H=\bar
H\comp H_\T$.
\end{enumerate}
\end{lem}
\begin{proof} 
The category $\T$ has weak kernels and therefore $\widehat\T$ is
abelian. Note that the weak kernel of a map $Y\to Z$ is obtained by
completing the map to an exact triangle $X\to Y\to Z\to \Si X$.  

Now let $H\colon\T\to\A$ be a cohomological functor. Extend $H$ to $\bar
H$ by sending $F$ in $\widehat\T$ with presentation
$$\T(-,X)\lto[(-,\p)] \T(-,Y)\lto F\lto 0$$ to the cokernel of $H\p$. The
functor $\bar H$ is automatically right exact, and it is exact because
$H$ is cohomological.
\end{proof}

The category $\widehat\T$ is called the {\em abelianization} of $\T$.
It has some special homological properties. Recall that an abelian
category is a {\em Frobenius category} if there are enough projectives
and enough injectives, and both coincide.
\begin{lem}
The abelianization of a triangulated category is an abelian Frobenius
category.
\end{lem}
\begin{proof}
The representable functors are projective objects in $\widehat\T$ by
Yoneda's lemma. Thus $\widehat\T$ has enough projectives.  Using the
fact that the Yoneda functors $\T\to\widehat\T$ and
$\T^\op\to\widehat{\T^\op}$ are universal cohomological functors, we
obtain an equivalence $\widehat\T^\op\to\widehat{\T^\op}$ which sends
$\T(-,X)$ to $\T(X,-)$ for all $X$ in $\T$. Thus the
representable functors are injective objects, and $\widehat\T$ has enough
injectives.
\end{proof}

\subsection{The idempotent completion of a triangulated category}
\label{ss:idpt}

Let $\T$ be a triangulated category. We identify $\T$ via the Yoneda
functor with a full subcategory of projective objects in the
abelianization $\widehat\T$. Recall that an additive category has {\em
split idempotents} if every idempotent map $\p^2=\p\colon X\to X$ has
a kernel. Note that in this case $X=\Ker\p\amalg\Ker (\id_X-\p)$.

\begin{lem}
The full subcategory $\bar\T$ of all projective objects in
$\widehat\T$ is a triangulated category with respect to the class of
triangles which are direct factors of exact triangles in $\T$. The
category $\bar\T$ has split idempotents, and every exact functor
$\T\to\U$ into a triangulated category with split idempotents extends
(up to a unique isomorphism) uniquely to an exact functor
$\bar\T\to\U$.
\end{lem}
\begin{proof}
The proof is straightforward, except for the verification of the
octahedral axiom which requires some work; see \cite{BS}.
\end{proof}

\begin{exm}
Let $\A$ be an abelian category. Then $\bfD^b(\A)$ has split
idempotents; see \cite{BS}.
\end{exm}

\subsection{Homotopy colimits}\label{ss:hocolim}

Let $\T$ be a triangulated category and suppose countable coproducts
exist in $\T$. Let $$X_0\lto[\p_0]X_1\lto[\p_1]X_2\lto[\p_2]\cdots$$ be a
sequence of maps in $\T$. A {\em homotopy colimit} of this sequence is
by definition an object $X$ which occurs in an exact triangle
$$\coprod_{i\geqslant 0}X_i\lto[(\id-\p_i)]\coprod_{i\geqslant 0}X_i\lto
X\lto\Si (\coprod_{i\geqslant 0}X_i).$$
Here, the $i$th component of the map $(\id-\p_i)$ is the composite
$$X_i\lto[\smatrix{\id\\-\p_i}]X_i\amalg
X_{i+1}\lto[\inc]\coprod_{i\geqslant 0}X_i.$$ Note that a homotopy
colimit is unique up to a (non-unique) isomorphism.

\begin{exm}
(1) Let $\p\colon X\to X$ be an idempotent map in $\T$, and denote by $X'$
a homotopy colimit of the sequence
$$X\lto[\p]X\lto[\p]X\lto[\p]\cdots.$$ The canonical map $X\to X'$ gives
rise to a split exact triangle $X''\to X\to X'\to\Si X''$, where
$X''=\Ker\p$. Thus $0\to X''\to X\to X'\to 0$ is a split exact
sequence and $X\cong X'\amalg X''$.

(2) Let $\A$ be an additive category and suppose countable coproducts
exist in $\A$. For a complex $X$ in $\A$ and $n\in\bbZ$, denote by
$\t_nX$ the truncation of $X$ such that $(\t_nX)^p=0$ for $p<n$ and
$(\t_nX)^p=X^p$ for $p\geqslant n$. For each $n\in\bbZ$, there is a
sequence of canonical maps
$$\t_nX\lto \t_{n-1}X\lto \t_{n-2}X \lto\cdots$$ which are compatible
with the canonical maps $\p_p\colon\t_pX\to X$. If $X'$ denotes a
homotopy colimit, then the family $(\p_{n-i})_{i\geqslant 0}$ induces
an isomorphism $X'\to X$ in $\bfK(\A)$.
\end{exm} 

\subsection{Brown representability}\label{ss:brown}

Let $\T$ be a triangulated category with arbitrary coproducts.  An
object $S$ in $\T$ satisfying the following conditions is called a
{\em perfect generator}.
\begin{enumerate}
\item[(PG1)] There is no proper full triangulated subcategory of $\T$ which
contains $S$ and is closed under taking coproducts.
\item[(PG2)] Given a countable family of maps $X_i\to Y_i$ in $\T$
such that the map $\Hom_\T(S,X_i)\rightarrow\Hom_\T(S,Y_i)$ is
surjective for all $i$, the induced map
$$\Hom_\T(S,\coprod_iX_i)\longrightarrow\Hom_\T(S,\coprod_iY_i)$$ is
surjective.
\end{enumerate}

We have the following {\em Brown representability theorem} for
triangulated categories with a perfect generator. 

\begin{thm}[Brown representability]
Let $\T$ be a triangulated category with arbitrary coproducts and a
perfect generator.
\begin{enumerate}
\item A functor $F\colon\T^\op\to\Ab$ is cohomological and sends all
coproducts in $\T$ to products if and only if $F\cong\Hom_\T(-,X)$ for
some object $X$ in $\T$.
\item An exact functor $\T\to\U$ between triangulated categories preserves 
all coproducts if and only if it has a right adjoint.
\end{enumerate}
\end{thm}

We give a complete proof. This is based on the following lemma, which
explains the condition (PG2) and is independent from the triangulated
structure of $\T$.

\begin{lem}
Let $\T$ be an additive category with arbitrary coproducts and weak
kernels.  Let $S$ be an object in $\T$, and denote by $\S$ the full
subcategory of all coproducts of copies of $S$.
\begin{enumerate}
\item The category $\S$ has weak kernels and $\widehat{\S}$ is an
abelian category.
\item The assignment $F\mapsto F|_{\S}$ induces 
an exact functor $\widehat\T\lto\widehat{\S}$.
\item The functor $\T\to\widehat\S$ sending $X$ to $\T(-,X)|_\S$
preserves countable coproducts if and only if {\rm (PG2)} holds.
\end{enumerate}
\end{lem}
\begin{proof}
First observe that for every  $X$ in $\T$, there exists 
an {\em approximation} $X'\to X$ such that
$X'\in\S$ and $\T(T,X')\to \T(T,X)$ is surjective for all $T\in\S$.
Take $X'=\coprod_{\a\in \T(S,X)}S$ and the canonical map $X'\to X$.

(1) To prove that $\widehat{\S}$ is abelian, it is sufficient to show
that every map in $\S$ has a weak kernel; see (\ref{ss:coh}).  To
obtain a weak kernel of a map $Y\to Z$ in $\S$, take the composite of
a weak kernel $X\to Y$ in $\T$ and an approximation $X'\to X$.

(2) We need to check that for $F$ in $\widehat \T$, the restriction $F|_{\S}$
belongs to $\widehat{\S}$. It is sufficient to prove this for $F=\T(-,Y)$.
To obtain a presentation, let $X\to Y'$ be a weak kernel of an approximation 
$Y'\to Y$. The composite $X'\to Y'$ with an approximation $X'\to X$ gives an
exact sequence 
$$\T(-,X')|_{\S}\lto \T(-,Y')|_{\S}\lto F|_{\S}\lto 0,$$
which is a
presentation $$\S(-,X')\lto \S(-,Y')\lto F|_{\S}\lto 0$$
of $F|_\S$.
Clearly, the assignment $F\mapsto F|_{\S}$ is exact.

(3) We denote by $I\colon \S\to\T$ the inclusion and write $I_*\colon
\widehat\T\to\widehat\S$ for the restriction functor. Note that $I$
induces a right exact functor $I^*\colon \widehat\S\to\widehat\T$ by
sending each representable functor $\S(-,X)$ to $\T(-,X)$.

The functor $\T\to\widehat\S$ preserves countable coproducts if and
only if $I_*$ preserves countable coproducts. We have $I_*\comp
I^*\cong\Id_{\widehat S}$ and this implies that $I_*$ preserves
countable coproducts if and only if $\Ker I_*$ is closed under
countable coproducts, since $I_*$ induces an equivalence
$\widehat\T/\Ker I_*\to\widehat\S$; see \cite[III.2]{G}. Now observe
that $\Ker I_*$ being closed under countable coproducts is a
reformulation of the condition (PG2).
\end{proof}

\begin{proof}[Proof of the Brown representability theorem]
  (1) Fix a perfect generator $S$. First observe that we may assume
  $\Si S\cong S$. Otherwise, replace $S$ by $\coprod_{n\in\bbZ} \Si^n
  S$. This does not affect the condition (PG2). Also, it does not
  affect the condition (PG1), because a triangulated subcategory
  closed under countable coproducts is closed under direct factors;
  see (\ref{ss:hocolim}).

We construct inductively a sequence
$$X_0 \lto[\p_0] X_1 \lto[\p_1] X_2 \lto[\p_2] \cdots$$ of maps in
$\T$ and elements $\pi_i$ in $FX_i$ as follows. Let $X_0=0$ and
$\pi_0=0$. Let $X_1=S^{[FS]}$ be the coproduct of copies of $S$
indexed by the elements in $FS$, and let $\pi_1$ be the element
corresponding to $\id_{FS}$ in $FX_1\cong (FS)^{FS}$. Suppose we have
already constructed $\p_{i-1}$ and $\pi_i$ for some $i>0$. Let
$$K_i=\{\a\in \T(S,X_i)\mid (F\a)\pi_{i}=0\}$$ and complete
the canonical map $\chi_i\colon  S^{[K_i]}\to X_i$ to an exact triangle
$$S^{[K_i]}\lto[\chi_i] X_{i} \lto[\p_{i}] X_{i+1} \lto\Si
S^{[K_i]}.$$ Now choose an element $\pi_{i+1}$ in $FX_{i+1}$ such that
$(F\p_{i})\pi_{i+1}=\pi_i$. This is possible since
$(F\chi_{i})\pi_{i}=0$ and $F$ is cohomological.

We identify each $\pi_i$ via Yoneda's lemma with a map $\T(-,X_i)\to
F$ and obtain the following commutative diagram with exact rows in
$\widehat\S$, where $\S$ denotes the full subcategory of all
coproducts of copies of $S$ in $\T$ and $\psi_i=\T(-,\p_i)|_\S$.
$$\xymatrix{ 0\ar[r]&\Ker\pi_i|_\S\ar[d]^0\ar[r]
&\T(-,X_i)|_\S\ar[d]^{\psi_i}\ar[r]^-{\pi_i} &F|_\S\ar@{=}[d]\ar[r]& 0\\
0\ar[r]&\Ker\pi_{i+1}|_\S\ar[r] &\T(-,X_{i+1})|_\S\ar[r]^-{\pi_{i+1}}
&F|_\S\ar[r]& 0}$$ We compute the colimit of the sequence
$(\psi_i)_{i\geqslant 0}$ in $\widehat\S$ and obtain an exact sequence
\begin{equation}\label{eq:brown}
0\lto\coprod_i\T(-,X_i)|_\S\lto[(\id-\psi_i)]\coprod_i\T(-,X_i)|_\S\lto
F|_\S\lto 0
\end{equation} 
because the sequence $(\psi_i)_{i\geqslant 0}$ is a coproduct of a
sequence of zero maps and a sequence of identity maps.

Next consider the exact triangle
$$\coprod_i X_i\lto[(\id-\p_i)]\coprod_i X_i\lto X\lto\Si(\coprod_i
X_i)$$ and observe that $$(\pi_i)\in\prod_i FX_i\cong F(\coprod_i X_i)$$
induces a map
$$\pi\colon \T(-,X)\lto F$$ by Yoneda's lemma. 
We have an isomorphism
$$\coprod_i\T(-,X_i)|_\S\cong\T(-,\coprod_iX_i)|_\S$$
because of the
reformulation of condition (PG2) in lemma~(\ref{ss:brown}), and we
obtain in $\widehat\S$ the following exact sequence.
\begin{multline*}
\coprod_i\T(-, X_i)|_\S\xto{(\id-\psi_i)}\coprod_i\T(-, X_i)|_\S\lto
\T(-,X)|_\S\lto \\ \coprod_i \T(-,\Si
X_i)|_\S\xto{(\id-\Si\psi_i)}\coprod_i\T(-,\Si X_i)|_\S 
\lto\T(-,\Si X)|_\S
\end{multline*}
A comparison with the exact sequence (\ref{eq:brown}) shows that
$$\pi|_\S\colon \T(-,X)|_\S\lto F|_\S$$ is an isomorphism, because
$(\id-\Si\psi_i)$ is a monomorphism. Here, we use that $\Si S\cong S$.

Next observe that the class of objects $Y$ in $\T$ such that $\pi_Y$
is an isomorphism forms a triangulated subcategory of $\T$ which is
closed under taking coproducts. Using condition (PG1), we conclude
that $\pi$ is an isomorphism.

(2) Let $F\colon\T\to\U$ be an exact functor. If $F$ preserves all
coproducts, then one defines the right adjoint $G$ by sending an object
$X$ in $\U$ to the object in $\T$ representing $\Hom_\U(F-,X)$.
Thus $$\Hom_\U(F-,X)\cong\Hom_\T(-,GX).$$
Conversely, given a right adjoint of $F$, it is automatic that $F$
preserves all coproducts.
\end{proof}

\begin{rem} 
(1) In the presence of (PG2), the condition (PG1) is equivalent to the
    following condition.
\begin{enumerate}
\item[(PG1')] Let $X$ be in $\T$ and suppose $\Hom_\T(\Si^nS,X)=0$ for
all $n\in\bbZ$. Then $X=0$.
\end{enumerate}

(2) The Brown representability theorem implies that a triangulated
category $\T$ with a perfect generator has arbitrary products. In
fact, given a family of objects $X_i$ in $\T$, let $\prod_iX_i$ be the
object representing $\prod_i\Hom_\T(-,X_i)$.

(3) There is the dual concept of a {\em perfect cogenerator} for a
    triangulated category. The dual Brown representability theorem for
    triangulated categories $\T$ with a perfect cogenerator
    characterizes the representable functors $\Hom_\T(X,-)$ as the
    cohomological and product preserving functors $\T\to\Ab$.
\end{rem}

\subsection{Notes}
The abelianization of a triangulated category appears in Verdier's
th\`ese \cite{V} and in Freyd's work on the stable homotopy category
\cite{F}. Note that their construction is slightly different from the
one given here, which is based on coherent functors in the sense of
Auslander \cite{A}. The idempotent completion of a triangulated
category is studied by Balmer and Schlichting in \cite{BS}.  Homotopy
colimits appear in work of B\"okstedt and Neeman \cite{BN}. The Brown
representability theorem in homotopy theory is due to Brown \cite{B}.
Generalizations of the Brown representability theorem for triangulated
categories can be found in work of Franke \cite{Fr}, Keller \cite{Ke2},
and Neeman \cite{N4,N3}. The proof given here follows \cite{Kr}.

\section{Resolutions}

Resolutions are used to replace a complex in some abelian category
$\A$ by another one which is quasi-isomorphic but easier to
handle. Depending on properties of $\A$, injective and projective
resolutions are constructed via Brown representability.

\subsection{Injective resolutions}\label{ss:ires}
Let $\A$ be an abelian category.  Suppose that $\A$ has arbitrary
products which are exact, that is, for every family of exact sequences
$X_i\to Y_i\to Z_i$ in $\A$, the sequence $\prod_i X_i\to \prod_i
Y_i\to\prod_i Z_i$ is exact. Suppose in addition that $\A$ has an
injective cogenerator which we denote by $U$. Thus $\Hom_\A(X,U)=0$
implies $X=0$ for every object $X$ in $\A$. Observe that
\begin{eqnarray}\label{eq:H}
\Hom_{\bfK(\A)}(-,U)\cong\Hom_\A(H^0-,U).
\end{eqnarray} 

Denote by $\bfK_\inj(\A)$ the smallest full triangulated subcategory
of $\bfK(\A)$ which is closed under taking products and contains all
injective objects of $\A$ (viewed as complexes concentrated in degree
zero).

\begin{lem}
The triangulated category $\bfK_\inj(\A)$ is perfectly cogenerated by
$U$.  Therefore the inclusion $\bfK_\inj(\A)\to\bfK(\A)$ has a left
adjoint $\bfi\colon\bfK(\A)\to\bfK_\inj(\A)$.
\end{lem}
\begin{proof}
Fix a family of maps $X_i\to Y_i$ in $\bfK(\A)$ and suppose that
$$\Hom_{\bfK(\A)}(Y_i,U)\lto\Hom_{\bfK(\A)}(X_i,U)$$
is surjective for
all $i$.  The isomorphism (\ref{eq:H}) shows that $H^0X_i\to H^0Y_i$
is a monomorphism for all $i$. Taking products in $\A$ is exact, and
therefore $H^0\colon\bfK(\A)\to\A$ preserves products. Thus
$H^0(\prod_i X_i)\to H^0(\prod_i Y_i)$ is a monomorphism and
$$\Hom_{\bfK(\A)}(\prod_i Y_i,U)\lto\Hom_{\bfK(\A)}(\prod_i X_i,U)$$
is surjective.

The existence of a left adjoint for the inclusion
$\bfK_\inj(\A)\to\bfK(\A)$ follows from the Brown representability
theorem
\end{proof}

The left adjoint $\bfi\colon \bfK(\A)\to\bfK_\inj(\A)$ induces for
each complex $X$ in $\A$ a natural map $X\to\bfi X$, and we may think
of this as an injective resolution. Recall that for $A$ in $\A$, a map
$A\to I$ of complexes in $\A$ is an {\em injective resolution} if it
is a quasi-isomorphism, each $I^n$ is injective, and $I^n=0$ for
$n<0$. An injective resolution $A\to I$ induces an isomorphism
$$\Hom_{\bfK(\A)}(I,Y)\cong\Hom_{\bfK(\A)}(A,Y)$$ for every complex
$Y$ in $\A$ with injective components. Moreover, $I$ belongs to
$\bfK_{\inj}(\A)$. Therefore $I\cong \bfi A$ in $\bfK(\A)$.

\begin{prop}
Let $\A$ be an abelian category.  Suppose $\A$ has an injective
cogenerator and arbitrary products which are exact. Let $X,Y$ be
complexes in $\A$.
\begin{enumerate}
\item The natural map $X\to\bfi X$ is a quasi-isomorphism and we have
natural isomorphisms 
\begin{eqnarray}\label{eq:i}
\Hom_{\bfD(\A)}(X,Y)\cong\Hom_{\bfD(\A)}(X,\bfi
Y)\cong\Hom_{\bfK(\A)}(X,\bfi Y).
\end{eqnarray}
\item The composite
$$\bfK_\inj(\A)\lto[\inc]\bfK(\A)\lto[\can]\bfD(\A)$$
is an equivalence of triangulated categories.
\end{enumerate}
\end{prop}

\begin{proof}
(1) The natural map $X\to \bfi X$ induces an isomorphism
$$\Hom_{\bfK(\A)}(\bfi X,U)\lto\Hom_{\bfK(\A)}(X,U)$$
and the
isomorphism (\ref{eq:H}) shows that $X\to\bfi X$ is a
quasi-isomorphism.

Now consider the class of complexes $Y'$ such that the map
$\Hom_{\bfK(\A)}(X,Y')\to\Hom_{\bfD(\A)}(X, Y')$ is bijective. This
class contains $U$, by lemma~(\ref{ss:ext}), and the objects form a
triangulated subcategory closed under taking products. Thus the map
$\Hom_{\bfK(\A)}(X,\bfi Y)\to\Hom_{\bfD(\A)}(X,\bfi Y)$ is bijective
for all $Y$, and
$$\Hom_{\bfD(\A)}(X,Y)\cong\Hom_{\bfD(\A)}(X,\bfi Y)$$
since $Y\to\bfi Y$ is an isomorphism in $\bfD(\A)$.

(2) The first part of the proof shows that the functor is fully
faithful and, up to isomorphism, surjective on objects. A
quasi-inverse $\bfD(\A)\to\bfK_\inj(\A)$ is induced by $\bfi$. This
follows from the universal property of the localization functor
$\bfK(\A)\to\bfD(\A)$, since $\bfi$ sends quasi-isomorphisms to
isomorphisms.
\end{proof}

\begin{rem}
(1) Given an arbitrary abelian category $\A$, the formation of the derived
category $\bfD(\A)$, via the localization with respect to all
quasi-isomorphisms, leads to a category where maps between two objects
not necessarily form a set. Thus appropriate assumptions on $\A$ are
needed. An equivalence $\bfK_\inj(\A)\to\bfD(\A)$ implies that
$\Hom_{\bfD(\A)}(X,Y)$ is actually a set for any pair of complexes
$X,Y$ in $\A$.

(2) The isomorphism (\ref{eq:i}) shows that the assignment $X\mapsto
\bfi X$ induces a right adjoint for the canonical functor
$\bfK(\A)\to \bfD(\A)$.

(3) Let $\Inj\A$ denote the full subcategory of injective objects in
    $\A$. Then
    $$\bfK^+(\Inj\A)\subseteq\bfK_\inj(\A)\subseteq\bfK(\Inj\A).$$ For
the first inclusion, see (\ref{ss:hocolim}). Note that
$\bfK_\inj(\A)=\bfK(\Inj\A)$ if every object in $\A$ has finite
injective dimension.
\end{rem}

\begin{exm} 
(1) An abelian category with a projective generator has exact products. 

(2) Let $\A$ be the category of modules over a ring $\La$.  Then $\La$
is a projective generator and $\Hom_{\bbZ}(\La^\op,\bbQ/\bbZ)$ is an
injective cogenerator for $\A$. Therefore $\A$ has exact products and
exact coproducts.

(3) The category of quasi-coherent sheaves on the projective line
$\bfP^1_k$ over a field $k$ does not have exact products.
\end{exm}

\subsection{Projective Resolutions}
The existence of injective resolutions turns into the existence of
projective resolutions if one passes from an abelian category to its
opposite category. Because of this symmetry, complexes in an abelian
category $\A$ admit projective resolutions provided that $\A$ has a
projective generator and arbitrary coproducts which are exact. Keep
these assumptions on $\A$, and denote by $\bfK_{\proj}(\A)$ the smallest
full triangulated subcategory of $\bfK(\A)$ which is closed under
taking coproducts and contains all projective objects of $\A$. Then
one obtains a right adjoint $\bfp\colon\bfK(\A)\to\bfK_{\proj}(\A)$ of
the inclusion. For every complex $X$, the natural map $\bfp X\to X$
has the dual properties of the injective resolution $X\to\bfi X$. The
precise formulation is left to the reader. Note that any module
category has a projective generator and exact coproducts.

\subsection{Derived functors}

Injective and projective resolutions are used to define derived
functors. Let $F\colon\A\to\B$ be a functor between abelian
categories.  The {\em right derived functor} $\bfR F\colon\bfD(\A)\to
\bfD(\B)$ of $F$ sends a complex $X$ to $F(\bfi X)$, and the {\em left
derived functor} $\bfL F\colon\bfD(\A)\to \bfD(\B)$ sends a complex
$X$ to $F(\bfp X)$.

\begin{exm}
Let $\La$ and $\Ga$ be a pair of rings and ${_\La B}_\Ga$ a
bimodule. Then the pair of adjoint functors
$$\xymatrix{\Mod\La\ar@<-1ex>[rrr]_{T=-\otimes_\La B}&&&
\Mod\Ga\ar@<-1ex>[lll]_{H=\Hom_\Ga(B,-)}}$$
induces a pair of adjoint functors
$$\xymatrix{\bfD(\Mod\La)\ar@<-1ex>[rrr]_{\bfL T=-\otimes^\bfL_\La
B}&&& \bfD(\Mod\Ga)\ar@<-1ex>[lll]_{\bfR H=\RHom_\Ga(B,-)}}.$$ This
follows from the fact that $\bfL T$ and $\bfR H$ are composed from
three pairs
$$\xymatrix{\bfD(\Mod\La)\ar@<-1ex>[r]_{\bfp}&
\bfK(\Mod\La)\ar@<-1ex>[l]_{\can}\ar@<-1ex>[rr]_{\bfK( T)}&&
\bfK(\Mod\Ga)\ar@<-1ex>[ll]_{\bfK(
H)}\ar@<-1ex>[r]_{\can}&\bfD(\Mod\Ga) \ar@<-1ex>[l]_{\bfi} }$$ of
adjoint functors.
\end{exm}

\subsection{Notes}
Injective and projective resolutions are needed to construct derived
functors. The first application of this formalism is Grothendieck's
duality theory \cite{Ha}. Resolutions of unbounded complexes are
established by Spaltenstein in \cite{S}, and also by B\"okstedt and
Neeman in \cite{BN}. The existence proof given here has the advantage
that it generalizes easily to other settings, for instance to
differential graded modules.

\section{Differential graded algebras and categories}

Differential graded algebras arise as complexes with an additional
multiplicative structure, and differential graded categories are
differential graded algebras with several objects. 
The concept generalizes that of an ordinary associative algebra, and
we study in a similar way categories of modules and derived categories
for such differential graded algebras and categories.

\subsection{Differential graded algebras and modules}
A {\em differential graded algebra} or {\em dg algebra} is a
$\bbZ$-graded associative algebra $$A=\coprod_{n\in\bbZ}A^n$$ over some
commutative ring $k$, together with a  differential $d\colon A\to
A$, that is, a homogeneous $k$-linear map of degree $+1$ satisfying
$d^2=0$ and the Leibniz rule
$$d(xy)=d(x)y+(-1)^nxd(y)\quad\textrm{for all}\quad x\in
A^n\quad\textrm{and}\quad y\in A.$$ A {\em dg $A$-module} is a
$\bbZ$-graded (right) $A$-module $X$, together with a differential $d\colon
X\to X$, that is, a homogeneous $k$-linear map of degree $+1$
satisfying $d^2=0$ and the Leibniz rule
$$d(xy)=d(x)y+(-1)^nxd(y)\quad\textrm{for all}\quad x\in
X^n\quad\textrm{and}\quad y\in A.$$ A morphism of dg $A$-modules is an
$A$-linear map which is homogeneous of degree $0$ and commutes with
the differential. We denote by $\bfC_\dg(A)$ the category of dg
$A$-modules.

\begin{exm} 
(1) An associative algebra $\La$ can be viewed as a dg algebra $A$ if
    one defines $A^0=\La$ and $A^n=0$ otherwise. In this case
    $\bfC_\dg(A)=\bfC(\Mod\La)$.

(2) Let $X,Y$ be complexes in some additive category $\C$. Define a
    new complex $\HOM_\C(X,Y)$ as follows. The $n$th component is
$$\prod_{p\in\mathbb Z}\Hom_\C(X^p,Y^{p+n})$$ 
and the  differential is given by 
$$d^n(\p^p)=d_Y\comp \p^p-(-1)^n\p^{p+1}\comp d_X.$$ Note that
$$H^n\HOM_\C(X,Y)\cong\Hom_{\bfK(\C)}(X,\Si^nY)$$ because $\Ker d^n$
identifies with $\Hom_{\bfC(\C)}(X,\Si^nY)$ and $\Im d^{n-1}$ with the
ideal of null-homotopic maps $X\to \Si^n Y$.  The composition of
graded maps yields a dg algebra structure for
$$\END_\C(X)=\HOM_\C(X,X)$$ and $\HOM_\C(X,Y)$ is a dg module over
$\END_\C(X)$.
\end{exm}

A map $\p\colon X\to Y$ of dg $A$-modules is {\em null-homotopic} if
there is a map $\r\colon X\to Y$ of graded $A$-modules which is
homogeneous of degree $-1$ such that $\p=d_Y\comp \r+\r\comp d_X$. The
null-homotopic maps form an ideal and the {\em homotopy category}
$\bfK_\dg(A)$ is the quotient of $\bfC_\dg(A)$ with respect to this
ideal. The homotopy category carries a triangulated structure which is
defined as before for the homotopy category $\bfK(\A)$ of an additive
category $\A$.

A map $X\to Y$ of dg $A$-modules is a {\em quasi-isomorphism} if it
induces isomorphisms $H^nX\to H^nY$ in each degree. The {\em derived
category} of  the dg algebra $A$ is 
by definition the localization
$$\bfD_\dg(A)=\bfK_\dg(A)[S^{-1}]$$ of $\bfK_\dg(A)$ with respect to
the class $S$ of all quasi-isomorphisms. Note that $S$ is a
multiplicative system and compatible with the
triangulation. Therefore $\bfD_\dg(A)$ is triangulated and the
localization functor $\bfK_\dg(A)\to\bfD_\dg(A)$ is exact.

\subsection{Differential graded categories}
Let $k$ be a commutative ring. A category is {\em $k$-linear} if the
maps between any two objects form a $k$-module and all composition
maps are bilinear. A category is {\em small} if the isomorphism
classes of objects form a set. Now fix a small $k$-linear category
$\A$. We think of $\A$ as an {\em algebra with several objects},
because a $k$-algebra is nothing but a $k$-linear category with
precisely one object. The {\em modules} over $\A$ are by definition
the $k$-linear functors $\A^\op\to\Mod k$. For example, the {\em free}
$\A$-modules are the representable functors
$$\A(-,A)=\Hom_\A(-,A)$$ with $A$ in $\A$. This terminology is
justified by the fact that the modules over an algebra $\La$ can be
identified with $k$-linear functors $\La^\op\to\Mod k$, where $\La$ is
viewed as a category with a single object.

A {\em dg category} $\A$ is by definition a dg algebra with several
objects. More precisely, $\A$ is a $\bbZ$-graded $k$-linear category,
that is,
$$\A(A,B)=\coprod_{n\in\bbZ}\A(A,B)^n$$ is a $\bbZ$-graded $k$-module
for all $A,B$ in $\A$, and the composition maps
$$\A(A,B)\times\A(B,C)\lto\A(A,C)$$ are bilinear and homogeneous of
degree $0$. In addition, there are differentials $d\colon \A(A,B)\to
\A(A,B)$ for all $A,B$ in $\A$, that is, homogeneous $k$-linear maps
of degree $+1$ satisfying $d^2=0$ and the Leibniz rule. 

The {\em opposite category} $\A^\op$ of a dg category $\A$ has the
same objects, the maps are $\A^\op(A,B)=\A(B,A)$ for all objects
$A,B$, and the composition maps are
$$\A^\op(A,B)^p\times\A^\op(B,C)^q\lto\A^\op(A,C)^{p+q}, \quad
(\phi,\psi)\mapsto (-1)^{pq}\phi\comp\psi.$$ Keeping the differentials
from $\A$, one checks easily that the Leibniz rule holds in $\A^\op$.

Let $\A$ be a small dg category. A {\em dg $\A$-module} is a graded
functor $X\colon \A^\op\to\Gr k$ into the category of graded $k$-modules, that
is, the maps $$\A^\op(A,B)\lto \Hom_{\Gr k}(XA,XB)$$ are homogeneous
$k$-linear of degree $0$, where $$\Hom_{\Gr
k}(M,N)^p=\{\p\in\Hom_k(M,N)\mid\p(M^q)\subseteq N^{q+p} \textrm{ for
all } q\in\bbZ\}$$ for $M,N$ in $\Gr k$. In addition, there are
differentials $d\colon XA\to XA$ for all $A$ in $\A$, that is,
homogeneous $k$-linear maps of degree $+1$ satisfying $d^2=0$ and the
Leibniz rule
$$d(xy)=d(x)y+(-1)^nxd(y)\quad\textrm{for all}\quad x\in
(XA)^n\quad\textrm{and}\quad y\in \A(B,A)^p,$$
where $xy=(-1)^{np}(Xy)x$.

\begin{exm}
(1) A dg algebra can be viewed as a dg category with one
    object. Conversely, a dg category $\A$ with one object $A$ gives a
    dg algebra $\A(A,A)$.

(2) A class $\A_0$ of complexes in some additive category $\C$ gives
rise to a dg category $\A$. Define $\A(X,Y)=\HOM_\C(X,Y)$ for $X,Y$ in
$\A_0$.  The composition in $\A$ is induced by the composition of graded
maps in $\C$.  If $\A$ is small and $Y$ is any complex in $\C$, the
functor sending $X$ in $\A$ to $\HOM_\C(X,Y)$ is a dg $\A$-module.
\end{exm}

Let $\A$ be a small dg category. A map $\p\colon X\to Y$ between dg
$\A$-modules is a natural transformation such that the maps
$\p_A\colon XA\to YA$ are homogeneous of degree zero and commute with
the differentials.  We denote by $\bfC_\dg(\A)$ the category of dg
$\A$-modules. The homotopy category $\bfK_\dg(\A)$ and the derived
category $\bfD_\dg(\A)$ are defined as before for dg algebras.

\subsection{Duality}
Let $\A$ be a small dg category over some commutative ring $k$. We fix
an injective cogenerator $E$ of the category of $k$-modules. The
duality $\Hom_k(-,E)$ between $k$-modules induces a duality between
dg modules over $\A$ and $\A^\op$ as follows. Let $X$ be a dg
$\A$-module. Define a dg $\A^\op$-module $DX$ by
$$\big((DX)A\big)^n=\Hom_k\big((XA)^{-n},E\big)\quad\textrm{and}\quad
\big((DX)\p\big)\a=(-1)^{pq}\a\comp X\p$$
for $\p\in\A(A,B)^p$ and $\a\in \big((DX)A\big)^q$.
The differential is $$d_{DX}^n=(-1)^{n+1}\Hom_k(d_X^{-n-1},E).$$  

Given $k$-modules $M,N$, there is a natural isomorphism
$$\Hom_k\big(M,\Hom_k(N,E)\big)\cong\Hom_k\big(N,\Hom_k(M,E)\big)$$ which sends
$\chi$ to $\Hom_k(\chi,E)\comp\d_N$. Here,
$$\d_N\colon N\lto\Hom_k\big(\Hom_k(N,E),E)\big)$$ is the evaluation map.
The isomorphism induces for dg modules $X$ over $\A$ and $Y$ over $\A^\op$ a
natural isomorphism
$$\Hom_{\bfC_\dg(\A)}(X,DY)\cong\Hom_{\bfC_\dg(\A^\op)}(Y,DX).$$
This isomorphism preserves null-homotopic maps and quasi-isomorphisms.

\begin{lem}
Let $X$ and $Y$ be dg modules over $\A$ and $\A^\op$
respectively. Then there are natural isomorphisms
\begin{align*}
\Hom_{\bfK_\dg(\A)}(X,DY)&\cong\Hom_{\bfK_\dg(\A^\op)}(Y,DX),\\
\Hom_{\bfD_\dg(\A)}(X,DY)&\cong\Hom_{\bfD_\dg(\A^\op)}(Y,DX).
\end{align*}
\end{lem}

\begin{rem}
The duality $D$ is unique up to isomorphism and depends only on $k$,
if we choose for $E$ a minimal injective cogenerator. Note that $E$
is {\em minimal} if it is an injective envelope of a coproduct of a
representative set of all simple $k$-modules.
\end{rem}

\subsection{Injective and projective resolutions}\label{ss:dgres}
Let $\A$ be a small dg category over some commutative ring $k$. For
$A,B$ in $\A$, we define the {\em free} dg $\A$-module
$A^\wedge=\A^\op(A,-)$ and the {\em injective} dg $\A$-module
$B^\vee=DB^\wedge$, where $B^\wedge$ denotes the free $\A^\op$-module
corresponding to $B$.

\begin{lem} 
Let $X$ be a dg $\A$-module. Then we have for $A\in\A$
$$\Hom_{\bfK_\dg(\A)}(A^\wedge,X)\cong H^0(XA)\quad\textrm{and}\quad
\Hom_{\bfK_\dg(\A)}(X,A^\vee)\cong \Hom_k(H^0(XA),E).$$
\end{lem}
\begin{proof}
The first isomorphism follows from Yoneda's lemma. In fact, taking
morphisms of graded $\A$-modules, we have $\Hom_\A(A^\wedge,X)\cong
XA$. Restricting to morphisms which commute with the differential, we
get $\Hom_{\bfC_\dg(\A)}(A^\wedge,X)\cong Z^0(XA)$. Finally, taking
morphism up to null-homotopic maps, we get the isomorphism
$\Hom_{\bfK_\dg(\A)}(A^\wedge,X)\cong H^0(XA)$.

The second isomorphism is an immediate consequence if one uses the
isomorphism
$$\Hom_{\bfK_\dg(\A)}(X,DA^\wedge)\cong\Hom_{\bfK_\dg(\A^\op)}(A^\wedge,DX).$$
\end{proof}

Let us denote by $\bfK_\pdg(\A)$ the smallest full triangulated
subcategory of $\bfK_\dg(\A)$ which is closed under coproducts and
contains all free dg modules $A^\wedge$, $A\in\A$.  Analogously,
$\bfK_\idg(\A)$ denotes the smallest full triangulated subcategory of
$\bfK_\dg(\A)$ which is closed under products and contains all
injective dg modules $A^\vee$, $A\in\A$.

\begin{lem}
\begin{enumerate} 
\item The category $\bfK_\pdg(\A)$ is perfectly generated by
$\coprod_{A\in\A}A^\wedge$. Therefore the inclusion has a right
adjoint $\bfp\colon\bfK_\dg(\A)\to\bfK_\pdg(\A)$.
\item The category $\bfK_\idg(\A)$ is perfectly cogenerated by
$\prod_{A\in\A}A^\vee$. Therefore the inclusion has a left adjoint
$\bfi\colon\bfK_\dg(\A)\to\bfK_\idg(\A)$.
\end{enumerate}
\end{lem}
\begin{proof} Adapt the proof of lemma~(\ref{ss:ires}).
\end{proof}

\begin{thm}
Let $\A$ be a small dg category and $X,Y$ be dg $\A$-modules.
\begin{enumerate}
\item The natural map $\bfp X\to X$ is a quasi-isomorphism and we have
$$\Hom_{\bfD_\dg(\A)}(X,Y)\cong\Hom_{\bfD_\dg(\A)}(\bfp
X,Y)\cong\Hom_{\bfK_\dg(\A)}(\bfp X, Y).$$
\item We have for all $A\in\A$ and $n\in\bbZ$
$$\Hom_{\bfD_\dg(\A)}(A^\wedge,\Si^nY)\cong H^n(YA).$$
\item The natural map $Y\to \bfi Y$ is a quasi-isomorphism and we have
$$\Hom_{\bfD_\dg(\A)}(X,Y)\cong\Hom_{\bfD_\dg(\A)}(X, \bfi
Y)\cong\Hom_{\bfK_\dg(\A)}(X, \bfi Y).$$
\item The functors
$$\bfK_\pdg(\A)\lto[\inc]\bfK_\dg(\A)\lto[\can]\bfD_\dg(\A)\quad\textrm{and}
\quad \bfK_\idg(\A)\lto[\inc]\bfK_\dg(\A)\lto[\can]\bfD_\dg(\A)$$ are
equivalences of triangulated categories.
\end{enumerate}
\end{thm}
\begin{proof} Adapt the proof of proposition~(\ref{ss:ires}).
\end{proof}

\subsection{Compact objects and perfect complexes}\label{ss:compact}

Let $\T$ be a triangulated category with arbitrary coproducts. An
object $X$ in $\T$ is called {\em compact}, if every map $X\to
\coprod_{i\in I} Y_i$ factors through $\coprod_{i\in J} Y_i$ for some
finite subset of $J\subseteq I$. Note that $X$ is compact if and only
if the functor $\Hom_\T(X,-)\colon\T\to\Ab$ preserves all
coproducts. This characterization implies that the compact objects in
$\T$ form a thick subcategory. Let us denote this subcategory by
$\T^c$.

\begin{lem} 
Suppose there is a set $\S$ of compact objects in $\T$ such that $\T$
admits no proper triangulated subcategory which contains $\S$ and is
closed under coproducts. Then $\T^c$ coincides with the smallest full
and thick subcategory which contains $\S$.
\end{lem}
\begin{proof}
First observe that the coproduct of all objects in $\S$ is a perfect
generator for $\T$.  Now let $X$ be a compact object in $\T$ and let
$F=\Hom_\T(-,X)$. Consider the sequence
$$X_0 \lto[\p_0] X_1 \lto[\p_1] X_2 \lto[\p_2] \cdots$$ of maps in the
proof of the Brown representability theorem. Thus $X$ is a homotopy
colimit of this sequence and fits into the exact triangle
$$\coprod_{i\geqslant 0}X_i\lto[(\id-\p_i)]\coprod_{i\geqslant
0}X_i\lto[\chi] X\lto[\psi]\Si (\coprod_{i\geqslant 0}X_i).$$ The map
$\psi$ factors through some finite coproduct $\coprod_i\Si X_i$ and
this implies $\psi=0$ because $(\id-\Si\p_i)\comp\psi=0$. Thus the
identity map $\id_X$ factors through $\chi$. Therefore $X$ is a direct
factor of some finite coproduct $\coprod_iX_i$. Each $X_i$ is obtained
from coproducts of objects in $\S$ by a finite number of extensions.
Comparing this with the construction of the thick subcategory
generated by $\S$ in (\ref{ss:tria}), one can show that $X$ belongs to
the thick subcategory generated by $\S$; see \cite[2.2]{N2} for
details.
\end{proof}

We can now describe the compact objects in some derived categories.
Let $\La$ be an associative ring.  Then $\La$, viewed as a complex
concentrated in degree zero, is a compact object in $\bfD(\Mod\La)$,
since $$\Hom_{\bfD(\Mod\La)}(\La,X)\cong H^0X.$$ A complex of
$\La$-modules is called {\em perfect} if it is quasi-isomorphic to a
bounded complex of finitely generated projective $\La$-modules.  The
perfect complexes are precisely the compact objects in $\bfD(\Mod\La)$.
This description of compact objects extends to the derived category
of a dg category.

\begin{prop}
Let $\A$ be a small dg category. A dg $\A$-module is compact in
$\bfD_\dg(\A)$ if and only if it belongs to the smallest thick
subcategory of $\bfD_\dg(\A)$ which contains all free dg modules
$A^\wedge$, $A\in\A$.
\end{prop}
\begin{proof} 
It follows from theorem~(\ref{ss:dgres}) that each free module is
  compact, and that $\bfD_\dg(\A)$ admits no proper triangulated
  subcategory which contains all free module and is closed under
  coproducts. Now apply lemma~(\ref{ss:compact}).
\end{proof}

\subsection{Notes}
Differential graded algebras were introduced by Cartan in order to
study the cohomology of Eilenberg-MacLane spaces \cite{C}.
Differential graded categories and their derived categories provide a
conceptual framework for tilting theory; they are studied
systematically by Keller in \cite{Ke2}. In particular, the existence
of projective and injective resolutions for dg modules is proved in
\cite{Ke2}; see also \cite{AH}. The analysis of compact objects is due
to Neeman \cite{N2}.

\section{Algebraic triangulated categories}

Algebraic triangulated categories are triangulated categories which
arise from algebraic constructions. There is a generic construction
which assigns to an exact Frobenius category its stable category.  On
the other hand, every algebraic triangulated category embeds into
the derived category of some dg category.

\subsection{Exact categories}

Let $\A$ be an exact category in the sense of Quillen \cite{Q}.  Thus
$\A$ is an additive category, together with a distinguished class of
sequences
\begin{equation*}\label{eq:ses}
0\lto X\lto[\a] Y\lto[\b] Z\lto 0
\end{equation*}
which are called {\em exact}.  The exact sequences satisfy a number of
axioms. In particular, the maps $\a$ and $\b$ in each exact sequence
as above form a {\em kernel-cokernel pair}, that is $\a$ is a kernel
of $\b$ and $\b$ is a cokernel of $\a$. A map in $\A$ which arises as
the kernel in some exact sequence is called {\em admissible mono}; a
map arising as a cokernel is called {\em admissible epi}. A full
subcategory $\B$ of $\A$ is {\em extension-closed} if every exact
sequence in $\A$ belongs to $\B$ provided its endterms belongs to
$\B$.

\begin{rem}
(1) Any abelian category is exact with respect to the class of
all short exact sequences. 

(2) Any full and extension-closed subcategory $\B$ of an exact
category $\A$ is exact with respect to the class of sequences which
are exact in $\A$. 

(3) Any small exact category arises, up to an exact equivalence, as a
    full and extension-closed subcategory of a module category; see
    (\ref{ss:difex}).
\end{rem}

\subsection{Frobenius categories}

Let $\A$ be an exact category. An object $P$ is {\em projective} if
the induced map $\Hom_\A(P,Y)\to\Hom_\A(P,Z)$ is surjective for every
admissible epi $Y\to Z$. Dually, an object $Q$ is {\em injective} if
the induced map $\Hom_\A(Y,Q)\to\Hom_\A(X,Q)$ is surjective for every
admissible mono $X\to Y$.  The category $\A$ has {\em enough
projectives} if every object $Z$ admits an admissible epi $Y\to Z$
with $Y$ projective.  And $\A$ has {\em enough injectives} if every
object $X$ admits an admissible mono $X\to Y$ with $Y$
injective. Finally, $\A$ is called a {\em Frobenius category}, if $\A$
has enough projectives and enough injectives and if both coincide.

\begin{exm}
(1) Let $\A$ be an additive category. Then $\A$ is an exact category
    with respect to the class of all split exact sequences in
    $\A$. All objects are projective and injective, and $\A$ is a
    Frobenius category.

(2) Let $\A$ be an additive category. The category $\bfC(\A)$ of
    complexes is exact with respect to the class of all sequences
    $0\to X\to Y\to Z\to 0$ such that $0\to X^n\to Y^n\to Z^n\to 0$ is
    split exact for all $n\in\bbZ$. A typical projective and injective
    object is a complex of the form
    $$I_A\colon\;\; \cdots \lto 0\lto A\lto[\id] A\lto 0\lto\cdots$$
    for some $A$ in $\A$. There is an obvious admissible mono
    $X\to\prod_{n\in\bbZ}\Si^{-n}I_{X^n}$ and also an admissible epi
    $\coprod_{n\in\bbZ}\Si^{-n-1}I_{X^n}\to X$. Also,
    $$\coprod_{n\in\bbZ}\Si^{-n}I_{X^n}\cong
    \prod_{n\in\bbZ}\Si^{-n}I_{X^n}.$$
    Thus $\bfC(\A)$ is a Frobenius
    category.
\end{exm}

\subsection{The derived category of an exact category}

Let $\A$ be an exact category.  A complex $X$ in $\A$ is called {\em
acyclic} if for each $n\in\bbZ$ there is an exact sequence
$$0\lto Z^n\lto[\a^n] X^{n}\lto[\b^n] Z^{n+1}\lto 0$$ in $\A$ such
that $d^n_X=\a^{n+1}\comp\b^n$.  A map $Y\to Z$ of complexes in $\A$
is a {\em quasi-isomorphism}, if it fits into an exact triangle $X\to
Y\to Z\to\Si X$ in $\bfK(\A)$ where $X$ is isomorphic to an acyclic
complex. Note that this definition coincides with our previous one, if
$\A$ is abelian.  The class $S$ of quasi-isomorphism in $\bfK(\A)$ is
a multiplicative system and compatible with the triangulation. The
derived category of the exact category $\A$ is by definition the
localization
$$\bfD(A)=\bfK(\A)[S^{-1}]$$ with respect to $S$. 

\begin{exm} 
Let $\A$ be an additive category and view it as an exact category with
respect to the class of all split exact sequences. Then $\bfD(\A)=\bfK(\A)$.
\end{exm}

\subsection{The stable category of a Frobenius category}
\label{ss:stab}

Let $\A$ be a Frobenius category. The {\em stable category}
$\bfS(\A)$ is by definition the quotient of $\A$ with respect to
the ideal $\I$ of maps which factor through an injective object.  Thus
$$\Hom_{\bfS(\A)}(X,Y)=\Hom_\A(X,Y)/\I(X,Y)$$ for all $X,Y$ in
$\A$. We choose for each $X$ in $\A$ an exact sequence
$$0\lto X\lto E\lto \Si X\lto 0$$ such that $E$ is injective, and
obtain an equivalence $\Si\colon\bfS(\A)\to\bfS(\A)$ by
sending $X$ to $\Si X$.  
Every exact sequence $0\to X\to Y\to Z\to 0$
fits into a commutative diagram with exact rows
$$\xymatrix{ 0\ar[r]&X
\ar[r]^\a\ar@{=}[d]&Y\ar[r]^\b\ar[d]&Z\ar[r]\ar[d]^\g&0\\
0\ar[r]&X\ar[r]&E\ar[r]&\Si X\ar[r]&0}$$ such that $E$ is injective.
A triangle in $\bfS(\A)$ is by definition {\em exact} if it
isomorphic to a sequence of maps $$ X\lto[\a] Y\lto[\b] Z\lto[\g]\Si
X$$ as above.

\begin{lem} 
The stable category of a Frobenius category is triangulated.
\end{lem}
\begin{proof} 
  It is easy to verify the axioms, once one observes that every map in
  $\bfS(\A)$ can be represented by an admissible mono in $\A$. Note
  that a homotopy cartesian square can be represented by a pull-back
  and push-out square. This gives a proof for (TR4').
\end{proof}

\begin{exm}
The category of complexes $\bfC(\A)$ of an additive category $\A$ is a
Frobenius category with respect to the degree-wise split exact
sequences. The maps factoring through an injective object are
precisely the null-homotopic maps. Thus the stable category of
$\bfC(\A)$ coincides with the homotopy category $\bfK(\A)$. Note that
the triangulated structures which have been defined via mapping cones
and via exact sequences in $\bfC(\A)$ coincide.
\end{exm}

\subsection{Algebraic triangulated categories}\label{ss:algebraic}

A triangulated category $\T$ is called {\em algebraic} if it 
satisfies the equivalent conditions of the following lemma.

\begin{lem} 
For a triangulated category $\T$, the following are equivalent.
\begin{enumerate}
\item There is an exact equivalence $\T\to\bfS(\A)$ for some Frobenius
category $\A$.
\item There is a fully faithful exact functor $\T\to\bfK(\B)$ for some
additive category $\B$.
\item There is a fully faithful exact functor $\T\to\bfS(\C)$ for some
Frobenius category $\C$.
\end{enumerate}
\end{lem}
\begin{proof}
(1) $\Rightarrow$ (2): Let $\U$ denote the full subcategory of all
objects in $\bfK(\A)$ which are acyclic complexes with injective
components. The functor sending a complex $X$ to $Z^0X$ induces an
equivalence $\U\to\bfS(\A)$.  Composing a quasi-inverse with with the
equivalence $\T\to\bfS(\A)$ and the inclusion $\U\to\bfK(\A)$ gives 
a fully faithful exact functor $\T\to\bfK(\A)$.

(2) $\Rightarrow$ (3): Clear, since $\bfK(\B)=\bfS(\bfC(\B))$.

(3) $\Rightarrow$ (1): We identify $\T$ with a full triangulated
subcategory of $\bfS(\C)$.  Denote by $F\colon \C\to\bfS(\C)$ the
canonical projection. Then the full subcategory $\A$ of all objects
$X$ with $FX$ in $\T$ is an extension-closed subcategory of $\C$,
containing all projectives and injectives from $\C$. Thus $\A$ is a
Frobenius category and $\bfS(\A)$ is equivalent to $\T$.
\end{proof}

\begin{exm} (1) Let $\A$ be an abelian category with enough injectives.
Then $\bfD^b(\A)$ is algebraic. To see this, 
observe that the composite
$$\bfK^+(\Inj\A)\lto[\inc]\bfK(\A)\lto[\can]\bfD(\A)$$
is fully
faithful; see lemma~(\ref{ss:ext}). The image of this functor contains
$\bfD^b(\A)$, because we can identify it with the smallest full
triangulated subcategory of $\bfD(\A)$ which contains the injective
resolutions of all objects in $\A$.  Thus $\bfD^b(\A)$ is equivalent
to a full triangulated subcategory of $\bfK(\A)$.

(2) The derived category $\bfD_\dg(\A)$ of a dg category $\A$ is
algebraic. This follows from the fact that the homotopy category
$\bfK_\dg(\A)$ is algebraic, and that $\bfD_\dg(\A)$ is equivalent to
the full triangulated subcategory $\bfK_\pdg(\A)$.
\end{exm}

\begin{thm}
  Let $\T$ be an algebraic triangulated category and $\S$ be a full
  subcategory which is small.  Then there exists a dg category
  $\A$ and an exact functor $F\colon \T\to \bfD_\dg(\A)$
  having the following properties.
\begin{enumerate}
\item $F$ identifies $\S$ with the full subcategory formed by the free
  dg modules $A^\wedge$, $A\in\A$.
\item $F$ identifies (up to direct factors) the smallest full and
  thick subcategory of $\T$ containing $\S$ with the full subcategory
  formed by the compact objects in $\bfD_\dg(\A)$.
\item Suppose $\T$ has arbitrary coproducts and every object in $\S$
  is compact. Then $F$ identifies the smallest full triangulated
  subcategory of $\T$ which contains $\S$ and is closed under all
  coproducts with $\bfD_\dg(\A)$.
\end{enumerate}
\end{thm}
\begin{proof}
Let $\T=\bfS(\B)$ for some Frobenius category $\B$. We denote by
$\tilde\B$ the full subcategory of $\bfC(\B)$ which is formed by all
acyclic complexes in $\B$ having injective components. This is a
Frobenius category with respect to the degree-wise split exact
sequences. The functor
$$\bfS(\tilde\B)\lto\bfS(\B),\quad X\mapsto Z^0X,$$
is an equivalence.

For each $X$ in $\T$, choose a complex $\tilde X$ in $\tilde\B$ with
$Z^0\tilde X=X$. Define $\A$ by taking as objects the set $\{\tilde
X\in\tilde\B\mid X\in\S\}$, and let $\A(\tilde X,\tilde
Y)=\HOM_\B(\tilde X,\tilde Y)$ for all $X,Y$ in $\S$. Then $\A$ is a
small dg category. We obtain a functor
$$\tilde\B\lto \bfC_\dg(\A),\quad X\mapsto \HOM_\B(-,X),$$ and compose it
with the canonical functor $\bfC_\dg(\A)\to\bfD_\dg(\A)$ to
get an exact functor
$$\bfS(\tilde\B)\lto \bfD_\dg(\A)$$ of triangulated categories. We
compose this with a quasi-inverse of $\bfS(\tilde\B)\to\bfS(\B)$ and
get the functor $$F\colon \T =\bfS(\B)\lto \bfD_\dg(\A), \quad
X\mapsto \HOM_\B(-,\tilde X).$$

(1) The functor $F$ sends $X$ in $\S$ to the free $\A$-module
    $\tilde X^\wedge$, and we have for $Y$ in $\S$
\begin{equation*}
\Hom_{\bfS(\B)}(X,Y)\cong H^0\HOM_{\B}(\tilde X,\tilde Y)\cong
H^0\A(\tilde X,\tilde Y)\cong\Hom_{\bfD_\dg(\A)}(\tilde
X^\wedge,\tilde Y^\wedge).
\end{equation*}
Thus the canonical map
\begin{equation}\label{eq:F}
\Hom_{\bfS(\B)}(X,Y)\lto\Hom_{\bfD_\dg(\A)}(FX,FY)
\end{equation}
is bijective for all $X,Y$ in $\S$.

(2) The functor $F$ is exact and therefore the bijection (\ref{eq:F})
for objects in $\S$ extends to a bijection for all objects in the
thick subcategory $\T'$ generated by $\S$.  Now observe that the
compact objects in $\bfD_\dg(\A)$ are precisely the objects in the
thick subcategory generated by the free dg $\A$-modules.  Thus $F$
identifies $\T'$, up to direct factors, with the full subcategory
formed by the compact objects in $\bfD_\dg(\A)$.

(3) The functor $F$ preserves coproducts. To see this, let $(X_i)$ be
a family of objects in $\T$ and $X$ in $\S$. The degree $n$ cohomology
of the canonical map
$$\coprod_i\HOM_\B(\tilde X,\tilde X_i)\lto\HOM_\B(\tilde
X,\coprod_i\tilde X_i)$$
identifies with the canonical map
$$\coprod_i\Hom_{\bfS(\B)}(X,\Si^nX_i)\lto
\Hom_{\bfS(\B)}(X,\Si^n\coprod_i X_i)$$ which is an isomorphism since
$X$ is compact. Thus $\coprod_iFX_i\cong F(\coprod_i X_i)$.  Using this
fact and the exactness of $F$, one sees that the map (\ref{eq:F}) is a
bijection for all $X,Y$ in the triangulated subcategory $\T''$ which
contains $\S$ and is closed under coproducts.  Now observe that
$\bfD_\dg(\A)$ is generated by all free dg $\A$-modules, that
is, there is no proper triangulated subcategory closed under
coproducts and containing all free dg $\A$-modules.  Thus $F$
identifies $\T''$ with $\bfD_\dg(\A)$.
\end{proof}

\begin{rem}
Let $\T$ be a small algebraic triangulated category and fix an
embedding $\T\to \bfD_\dg(\A)$. Suppose $S$ is a multiplicative system
of maps in $\T$ which is compatible with the triangulation.  Then the
localization functor $\T\to\T[S^{-1}]$ admits a fully faithful and
exact ``right adjoint'' $\T[S^{-1}]\to \bfD_\dg(\A)$; see
\cite{N2}. Therefore $\T[S^{-1}]$ is algebraic.
\end{rem}

\subsection{The stable homotopy category is not algebraic}

There are triangulated categories which are not algebraic. For
instance, the stable homotopy category of spectra is not algebraic.
The following argument has been suggested by Bill Dwyer. Given any
endomorphism $\p\colon X\to X$ in a triangulated category $\T$, denote
by $X/\p$ its cone. If $\T$ is algebraic, then we can identify $X$
with a complex and $\p$ induces an endomorphism of the mapping cone
$X/\p$  which is null-homotopic.  Thus $2\cdot \id_{X/\p}=0$ in
$\T$ for $\p=2\cdot\id_X$.  On the other hand, let $S$ denote the
sphere spectrum. Then it is well-known (and can be shown using
Steenrod operations) that the identity map of the mod $2$ Moore
spectrum $M(2)=S/(2\cdot \id_S)$ has order $4$.

\begin{prop}
There is no faithful exact functor from the stable homotopy category of spectra
into an algebraic triangulated category.
\end{prop}
\begin{proof} Let $F\colon\S\to\T$ be an exact functor
from the stable homotopy category of spectra to an algebraic
triangulated category and let $X=F(S)$. Then we have
$F(M(2))=X/(2\cdot \id_X)$ and therefore $F(2\cdot\id_{M(2)})=2\cdot
\id_{X/(2\cdot\id_X)}=0$. On the other hand, $2\cdot\id_{M(2)}\neq 0$. Thus $F$ is not
faithful.
\end{proof}

\subsection{The differential graded category of an exact category}
\label{ss:difex}

Let $\A$ be a small exact category. We denote by $\Lex\A$ the category
of additive functors $F\colon\A^\op\to\Ab$ to the category of abelian
groups which are {\em left exact}, that is, each exact sequence $0\to
X\to Y\to Z\to 0$ in $\A$ induces an exact sequence
$$0\lto FZ\lto FY\lto FX$$ of abelian groups.

\begin{lem} 
The category $\Lex\A$ is an abelian Grothendieck category. The
Yoneda functor
$$\A\lto\Lex\A,\quad X\mapsto\Hom_\A(-,X),$$ is exact and identifies $\A$
with a full extension-closed subcategory of $\Lex\A$.
It induces a fully faithful exact functor
$$\bfD^b(\A)\lto\bfD^b(\Lex\A).$$
\end{lem}
\begin{proof} 
For the first part, see \cite[A.2]{Ke1}.  We identify $\A$ with its
image in $\Lex\A$ under the Yoneda functor.  The proof of the first
part shows that for each exact sequence $0\to X\to Y\to Z\to 0$ in
$\Lex\A$ with $Z$ in $\A$, there exists an exact sequence $0\to X'\to
Y'\to Z\to 0$ in $\A$ such that the map $Y'\to Z$ factors through the
map $Y\to Z$. This implies that the category $\bfK^-(\A)$ is {\em left
cofinal} in $\bfK^-(\Lex\A)$ with respect to the class of
quasi-isomorphisms; see \cite[III.2.4.1]{V} or \cite[12.1]{Ke3}.  Thus
for every quasi-isomorphism $\s\colon X\to X'$ in $\bfK^-(\Lex\A)$
with $X'$ in $\bfK^-(\A)$, there exists a map $\p\colon X''\to X$ with
$X''$ in $\bfK^-(\A)$ such that $\s\comp\p$ is a quasi-isomorphism.
Therefore the canonical functor $\bfD^-(\A)\to\bfD^-(\Lex\A)$ is fully
faithful.
\end{proof}

An abelian Grothendieck category has enough injective objects.  Thus
$\bfD^b(\Lex\A)$ is algebraic, and therefore $\bfD^b(\A)$ is
algebraic.  This has the following consequence.

\begin{prop}
Let $\A$ be a small exact category. Then there exists a 
dg category $\bar\A$ and a fully faithful exact functor
$F\colon \bfD^b(\A)\to \bfD_\dg(\bar\A)$ having the following properties.
\begin{enumerate}
\item $F$ identifies $\A$ with the
full subcategory formed by the  free dg modules $A^\wedge$, $A\in\bar\A$.
\item $F$ identifies $\bfD^b(\A)$ (up to direct factors) with the full
subcategory formed by the compact objects in $\bfD_\dg(\bar\A)$.
\end{enumerate}
\end{prop}

\begin{exm} 
  Let $\La$ be a right noetherian ring and denote by $\A=\mod\La$ the
  category of finitely generated $\La$-modules. Then the functor
$$\Mod\La\lto\Lex\A,\quad X\mapsto\Hom_\La(-,X)|_\A,$$
is an equivalence; see \cite[II.4]{G}. Thus the exact embedding 
$\mod\La\to\Mod\La$ induces a fully faithful exact functor
$$\bfD^b(\mod\La)\lto\bfD^b(\Mod\La),$$ which identifies
$\bfD^b(\mod\La)$ with the full subcategory of complexes $X$ in
$\bfD^b(\Mod\La)$ such that $H^nX$ is finitely generated for all $n$
and $H^nX=0$ for almost all $n\in\bbZ$.
\end{exm}

\subsection{Notes}

Frobenius categories and their stable categories appear in the work of
Heller \cite{He}, and later in the work of Happel on derived
categories of finite dimensional algebras \cite{H}.  The derived
category of an exact category is introduced by Neeman in
\cite{N1}. The characterization of algebraic triangulated categories
via dg categories is due to Keller \cite{Ke2}.

\begin{appendix}
\section{The octahedral axiom}
Let $\T$ be a pre-trianglated category. In this appendix, it is shown
that the octahedral axiom (TR4) is equivalent to the axiom (TR4').
It is convenient to introduce a further axiom.
\begin{enumerate}
\item[(TR4'')] 
Every diagram
$$\xymatrix{ X\ar[r]\ar[d]^{}&Y\ar[d]\ar[r]&Z\ar[r]&\Si X\\
X'\ar[r]^{}&Y' }$$ consisting of a homotopy cartesian square with
differential $\d\colon Y'\to\Si X$ and an exact triangle can be
completed to a morphism
$$\xymatrix{
  X\ar[r]\ar[d]^{}&Y\ar[r]\ar[d]^{}&Z\ar[r]\ar@{=}[d]^{}&\Si
  X\ar[d]^{}\\ X'\ar[r]^{}&Y'\ar[r]^{}&Z\ar[r]^{}&\Si X' }$$ between
  exact triangles such that the composite $Y'\to Z\to\Si X$ equals $\d$.
\end{enumerate}

\begin{prop} 
  Given a pre-triangulated category, the axioms (TR4), (TR4'), and
 (TR4'') are all equivalent.
\end{prop}
\begin{proof}
(TR4') $\Rightarrow$ (TR4''):
Suppose there is a diagram
$$\xymatrix{
U\ar[r]^\a\ar[d]^{\p}&V\ar[d]^{\psi}\ar[r]^\b&W\ar[r]^\g&\Si U\\
X\ar[r]^{\k}&Y }$$ consisting of a homotopy cartesian square with
differential $\d\colon Y\to\Si U$ and an exact triangle.  Applying
(TR4'), we obtain a morphism
$$\xymatrix{
  U\ar[r]^\a\ar[d]^{\p}&V\ar[d]^{\psi'}\ar[r]^{\b'}&W'\ar@{=}[d]\ar[r]^{\g'}&\Si
  U\ar[d]^{\Si\p}\\
  X\ar[r]^{\k'}&Y'\ar[r]^{\la'}&W'\ar[r]^{\mu'}&\Si X }$$
between
exact triangles such that the left hand square is homotopy cartesian
with differential $\g'\comp\la'$. We apply (TR3) and
lemma~(\ref{ss:uni}) to obtain an isomorphism $(\id_U,\id_V,\s)$
between $(\a,\b,\g)$ and $(\a,\b',\g')$.  This yields the following
morphism
$$\xymatrix{
U\ar[r]^\a\ar[d]^{\p}&V\ar[d]^{\psi'}\ar[r]^{\b}&W\ar@{=}[d]\ar[r]^{\g}&\Si
U\ar[d]^{\Si\p}\\
X\ar[r]^{\k'}&Y'\ar[r]^{\s^{-1}\comp\la'}&W\ar[r]^{\mu'\comp\s}&\Si X
}$$ between exact triangles. Next we apply (TR3) and
lemma~(\ref{ss:uni}) to obtain an isomorphism $(\id_U,\id_{V\amalg
X},\t)$ between the triangles $(\smatrix{\a\\
\p},\smatrix{\psi&\k},\d)$ and $(\smatrix{\a\\
\p},\smatrix{\psi'&\k'},\g'\comp\la')$.  Let
$\la=\s^{-1}\comp\la'\comp\t$.  This yields the following morphism
$$\xymatrix{
U\ar[r]^\a\ar[d]^{\p}&V\ar[d]^{\psi}\ar[r]^{\b}&W\ar@{=}[d]\ar[r]^{\g}&\Si
U\ar[d]^{\Si\p}\\
X\ar[r]^{\k}&Y\ar[r]^{\la}&W\ar[r]^{\mu'\comp\s}&\Si X }$$
between exact triangles. Note that $\g\comp\la$ is a differential for
the left hand square. Thus (TR4'') holds.

(TR4'') $\Rightarrow$ (TR4):
Suppose there are exact triangles $(\a_1,\a_2,\a_3)$,
$(\b_1,\b_2,\b_3)$, and $(\g_1,\g_2,\g_3)$ with
$\g_1=\b_1\comp\a_1$. Use (TR1) and complete the map $Y\to U\amalg Z$
to an exact triangle
$$Y\xto{\smatrix{\a_2\\ \b_1}} U\amalg Z\xto{\smatrix{\d_1&-\g_2}}
V\lto[\e] \Si Y.$$ Then we apply (TR4'') twice to obtain a commutative diagram
$$\xymatrix{X\ar[r]^{\a_1}\ar@{=}[d]&Y\ar[r]^{\a_2}\ar[d]^{\b_1}&
U\ar[r]^{\a_3}\ar[d]^{\d_1}& \Si X\ar@{=}[d]\\
X\ar[r]^{\g_1}&Z\ar[r]^{\g_2}\ar[d]^{\b_2}&
V\ar[r]^{\g_3}\ar[d]^{\d_2}&\Si X\\ &W\ar@{=}[r]\ar[d]^{\b_3}&
W\ar[d]^{\d_3}\\ &\Si Y\ar[r]^{\Si\a_2}&\Si U }$$ such that the first
two rows and the two central columns are exact triangles.  Moreover,
the top central square is homotopy cartesian and the differential
satisfies $$-\Si\a_1\comp\g_3=\e=-\b_3\comp\d_2.$$ Thus (TR4) holds.

(TR4) $\Rightarrow$ (TR4''):
Suppose there is given a diagram
$$\xymatrix{X\ar[r]^\a\ar[d]^{\p_1}&Y\ar[d]^{\p_2}\ar[r]^\b&Z\ar[r]^\g&\Si
X\\ X'\ar[r]^{\a'}&Y'}$$ consisting of a homotopy cartesian square and
an exact triangle.  Thus we have an exact triangle
$$X\xto{\smatrix{\a \\ \p_1}} Y\amalg
X'\xto{\smatrix{\p_2&-\a'}}Y'\lto[\d] \Si X.$$ We apply (TR4) and
obtain the following commutative diagram
$$\xymatrix{X\ar[r]^-{\smatrix{\a \\ \p_1}}\ar@{=}[d]&Y\amalg
X'\ar[r]^-{\smatrix{\p_2&-\a'}}\ar[d]^-{\smatrix{\id&0}}&
Y'\ar[r]^-{\d}\ar[d]^-{\b'}& \Si X\ar@{=}[d]\\
X\ar[r]^-{\a}&Y\ar[r]^-{\b}\ar[d]^-{0}&
Z\ar[r]^-{\g}\ar[d]^-{\g'}&\Si X\ar[d]^-{\Si\smatrix{\a\\ \p_1}}\\
&\Si X'\ar@{=}[r]\ar[d]^-{\smatrix{0\\ \id}}& \Si X'\ar[d]^-{-\Si\a'}\ar[r]^-{\smatrix{0\\ \id}}&\Si(Y\amalg X')\\ &\Si
(Y\amalg X')\ar[r]^-{\Si\smatrix{\p_2&-\a'}}&\Si Y' }$$
such that $(\b',\g',-\Si\a')$ is an exact triangle.
This gives the following morphism
$$\xymatrix{
X\ar[r]^\a\ar[d]^{\p_1}&Y\ar[r]^\b\ar[d]^{\p_2}&Z\ar[r]^\g\ar@{=}[d]&\Si
X\ar[d]^{\Si\p_1}\\ X'\ar[r]^{\a'}&Y'\ar[r]^{\b'}&Z\ar[r]^{\g'}&\Si X'
}$$ of triangles where $\d=\g\comp\b'$ is the differential of the
homotopy cartesian square. Thus (TR4'') holds. In particular, (TR4')
holds and therefore the proof is complete.
\end{proof}
\end{appendix}


\begin{thebibliography}{99}
%
\bibitem{A}{\sc M. Auslander:} Coherent functors, In: Proc. Conf.
  Categorical Algebra (La Jolla, Calif., 1965), Springer-Verlag, New
  York (1966), 189--231.
%
\bibitem{AH}{\sc L. Avramov and S. Halperin:} Through the looking
glass: a dictionary between rational homotopy theory and local
algebra.  In: Algebra, algebraic topology and their interactions
(Stockholm, 1983), Springer Lecture Notes in Math. {\bf 1183} (1986),
1--27.
%
\bibitem{BS} {\sc P. Balmer and M. Schlichting:} Idempotent completion
of triangulated categories. J. Algebra {\bf 236} (2001), 819--834.
%
\bibitem{BN}{\sc M. B\"okstedt and A. Neeman:} Homotopy limits in
triangulated categories. Compositio Math. {\bf 86} (1993) 209--234.
%
\bibitem{B}{\sc E. H. Brown:} Cohomology theories. Annals of Math. {\bf
    75} (1962), 467--484.
%
\bibitem{C}{\sc H. Cartan:} Alg\`ebres d'Eilenberg-MacLane. Expos\'es
2 \`a 11, S\'em. H. Cartan, \'Ec. Normale Sup. (1954-1955),
S\'ecr\'etariat Math., Paris (1956).
%
\bibitem{Fr}{\sc J. Franke:} On the Brown representability theorem for
triangulated categories. Topology {\bf 40} (2001), 667--680.
%
\bibitem{F} {\sc P. Freyd:} Stable homotopy. In:
Proc. Conf. Categorical Algebra (La Jolla, Calif., 1965), Springer-Verlag, New
York (1966), 121--172.
%
\bibitem{G}{\sc P. Gabriel:} Des cat\'egories
ab\'eliennes. Bull. Soc. Math. France {\bf 90} (1962), 323--448.
%
\bibitem{GZ}{\sc P. Gabriel and M. Zisman:} Calculus of fractions and
homotopy theory.  Ergebnisse der Mathematik und ihrer Grenzgebiete
{\bf 35}, Springer-Verlag, New York (1967).
%
\bibitem{GM}{\sc S. I. Gelfand and Yu. I. Manin:} Methods of
homological algebra, Springer-Verlag (1996).
%
\bibitem{H}{\sc D. Happel:} On the derived category of a finite
dimensional algebra. Comment. Math. Helv. {\bf 62} (1987), 339--389.
%
\bibitem{Ha}{\sc R. Hartshorne:} Residues and Duality. Springer
  Lecture Notes in Math. {\bf 20} (1966).
%
\bibitem{He}{\sc A. Heller:} The loop space functor in homological
algebra. Trans. Amer. Math. Soc. {\bf 96} (1960), 382--394.
%
\bibitem{KS}{\sc M. Kashiwara and P. Schapira:} Categories and sheaves, Springer-Verlag (2005).
%
\bibitem{Ke1}{\sc B. Keller:} Chain complexes and stable categories. 
Manus. Math. {\bf 67} (1990), 379--417.
%
\bibitem{Ke2}{\sc B. Keller:} Deriving DG categories.
Ann. Sci. \'Ecole. Norm. Sup. {\bf 27} (1994), 63--102.
%
\bibitem{Ke3}{\sc B. Keller:} Derived categories and their uses.
In: Handbook of algebra, Vol. 1, North-Holland (1996), 671--701. 
%
\bibitem{Kr}{\sc H. Krause:}
A Brown representability theorem via coherent functors. 
Topology {\bf 41} (2002), 853--861.
%
\bibitem{M}{\sc J. P. May:}
The additivity of traces in triangulated categories.
Adv. Math. {\bf 163} (2001), 34--73.
%
\bibitem{N1}{\sc A. Neeman:} The derived category of an exact category.
J. Algebra {\bf 135} (1990), 388--394.
%
\bibitem{N2}{\sc A. Neeman:} The connection between the $K$-theory
localization theorem of Thomason, Trobaugh and Yao and the smashing
subcategories of Bousfield and Ravenel. Ann. Sci. \'Ecole
Norm. Sup. {\bf 25} (1992), 547--566.
%
\bibitem{N4}{\sc A. Neeman:}
The Grothendieck duality theorem via Bousfield's techniques and 
Brown representability. J. Amer. Math. Soc. {\bf 9} (1996), 205--236.
%
\bibitem{N3}{\sc A. Neeman:} Triangulated categories. Annals of
Mathematics Studies {\bf 148}, Princeton University Press (2001).
%
\bibitem{PS}{\sc B. J. Parshall and L. L. Scott:} Derived categories,
quasi-hereditary algebras, and algebraic groups. Carleton
U. Math. Notes {\bf 3} (1988), 1--144.
%
\bibitem{P}{\sc D. Puppe:} On the structure of stable homotopy theory.
In:  Colloquium on algebraic topology. Aarhus Universitet Matematisk
  Institut (1962), 65--71.
%
\bibitem{Q}{\sc D. Quillen:} Higher algebraic K-theory, I.  In:
Algebraic K-theory, Springer Lecture Notes in Math. {\bf 341} (1973),
85--147.
%
\bibitem{S}{\sc N. Spaltenstein:} Resolutions of unbounded
complexes. Compositio Math. {\bf 65} (1988) 121--154.
%
\bibitem{V}{\sc J. L. Verdier:} Des cat\'egories d\'eriv\'ees des
cat\'egories ab\'eliennes. Ast\'erisque {\bf 239} (1996).
%
\bibitem{W}{\sc C. Weibel:} An introduction to homological algebra.
  Cambridge studies in advanced mathematics {\bf 38}, Cambridge
  University Press (1994).
%
\end{thebibliography}
\end{document}